\newif\iffinal
\finaltrue	

\documentclass[letterpaper,12pt,reqno]{amsart}
\RequirePackage[utf8]{inputenc}
\usepackage[portrait,margin=2.6cm]{geometry}
\iffinal\else\usepackage[notref,notcite]{showkeys}\fi
\usepackage{mathrsfs,soul}
\usepackage{hyperref}
\usepackage[foot]{amsaddr}
\usepackage{amssymb,amsthm,amsfonts,amsbsy,latexsym,dsfont}
\usepackage[textsize=small]{todonotes}       
\usepackage{graphicx}
\usepackage[numeric,initials,nobysame]{amsrefs}
\usepackage{upref,setspace}

\usepackage{caption}
\usepackage{subcaption}
\usepackage{fourier}

\usepackage{enumerate}

\numberwithin{equation}{section}
\numberwithin{figure}{section}
\numberwithin{table}{section}

 \let\go=w

\newcommand{\cC}{\mathcal{C}}
\newcommand{\cE}{\mathcal{E}}\newcommand{\cF}{\mathcal{F}}

\newcommand{\vzero}{\mathbf{0}}

\newcommand{\ve}{\mathbf{e}}

\newcommand{\bN}{\mathbb{N}}

\newcommand{\bZ}{\mathbb{Z}}

\newtheorem{thm}{Theorem}[section]
\newtheorem{lem}[thm]{Lemma}
\newtheorem{cor}[thm]{Corollary}

\newtheorem*{ass*}{Assumption}
\newtheorem*{theorem*}{Theorem}

\newtheorem{conj}[thm]{Conjecture}

\theoremstyle{definition}

\newtheorem{rem}{Remark}
\renewcommand{\leq}{\leqslant}
\renewcommand{\geq}{\geqslant}

\newcommand{\ind}{\mathds{1}}

\newcommand{\set}[1]{\left\{#1\right\}}

\newcommand{\convp}{\stackrel{\mathrm{P}}{\longrightarrow}}
\newcommand{\convd}{\stackrel{\mathrm{d}}{\Longrightarrow}}

\newcommand{\convas}{\stackrel{\mathrm{a.s.}}{\longrightarrow}}

\def\qed{ \hfill $\blacksquare$}
\DeclareMathOperator{\E}{\mathbb{E}}
\DeclareMathOperator{\pr}{\mathbb{P}}

\DeclareMathOperator{\var}{Var}
\DeclareMathOperator{\cov}{Cov}

\newcommand{\bpare}[1]{\left( #1 \right)}
\newcommand{\bbrac}[1]{\left[ #1 \right]}
\DeclareMathOperator{\Ann}{Ann}

\begin{document}

\title[Critical First Passage Percolation]{Asymptotics for $2D$ Critical First Passage Percolation}

\date{}

\author[Damron]{Michael Damron$^1$}
\address{$^1$Indiana University}
\author[Lam]{Wai-Kit Lam$^2$}
\address{$^2$Indiana University}
\author[Wang]{Xuan Wang$^3$}
\address{$^3$Indiana University and ICERM}
\email{mdamron6@gatech.edu, lamw@umail.iu.edu, xuanwang9527@gmail.com}

\maketitle

\begin{abstract}
We consider first-passage percolation on $\bZ^2$ with i.i.d. weights, whose distribution function satisfies $F(0) = p_c = 1/2$. This is sometimes known as the ``critical case'' because large clusters of zero-weight edges force passage times to grow at most logarithmically, giving zero time constant. Denote $T(\vzero, \partial B(n))$ as the passage time from the origin to the boundary of the box $[-n,n] \times [-n,n]$. We characterize the limit behavior of $T(\vzero, \partial B(n))$ by conditions on the distribution function $F$. We also give exact conditions under which $T(\vzero, \partial B(n))$ will have uniformly bounded mean or variance. These results answer several questions of Kesten and Zhang from the '90s and, in particular, disprove a conjecture of Zhang (\cite{zhang1999double}) from '99. In the case when both the mean and the variance go to infinity as $n \to \infty$, we prove a CLT under a minimal moment assumption. The main tool involves a new relation between first-passage percolation and invasion percolation: up to a constant factor, the passage time in critical first-passage percolation has the same first-order behavior as the passage time of an optimal path constrained to lie in an embedded invasion cluster.
\end{abstract}


\section{Introduction}
\subsection{The model}
Consider the integer lattice $\bZ^d$ and denote by $\mathcal{E}^d$ the set of nearest-neighbor edges. Given a distribution function $F$ with $F(0^-)=0$, let $(t_e : e \in \mathcal{E}^d)$ be a family of i.i.d. random variables (edge-weights) with common distribution function $F$. In first-passage percolation, we study the random pseudo-metric on $\mathbb{Z}^d$ induced by these edge-weights.

The model is defined as follows. For $x,y\in\bZ^d$, a (vertex self-avoiding) path from $x$ to $y$ is an alternating sequence $(v_0,e_1,v_1,\ldots,e_n,v_n)$, where the $v_i$'s, $i=1,\ldots,n-1$, are distinct vertices in $\bZ^d$ which are different from $x$ or $y$, and $v_0=x$, $v_n=y$; $e_i$ is an edge in $\mathcal{E}^d$ which connects $v_{i-1}$ and $v_i$. If $x=y$, the path is called a (vertex self-avoiding) circuit. For a path $\gamma$, we define the passage time of $\gamma$ to be $T(\gamma) = \sum_{e\in\gamma} t_e$. For any $A$, $B\subset \bZ^d$, we define the \emph{first-passage time} from $A$ to $B$ by
\[
T(A,B) = \inf\{T(\gamma): \gamma\text{ is a path from a vertex in } A \text{ to a vertex in }B\}.
\]
For $A=\{x\}$, write $T(x,B)$ for $T(\{x\},B)$ and similarly for $B$. A geodesic is a path $\gamma$ from $A$ to $B$ such that $T(\gamma) = T(A,B)$.

It is a consequence of the sub-additive ergodic theorem that if $\E T(x,y)<\infty$ for all $x,y$ then there exists a constant $\mu$, called the \emph{time constant}, such that
\[
\lim_{n\to\infty} \frac{T(\vzero,n\ve_1)}{n} = \mu \text{ almost surely and in }L^1\ ,
\]
where $\ve_1 = (1,0,\ldots, 0)$. It was shown by Kesten \cite[Theorem~6.1]{kesten1986aspects} that 
\begin{equation}\label{eq: old_kesten}
\mu=0 \text{ if and only if } F(0)\geq p_c,
\end{equation} 
where $p_c$ is the critical probability for Bernoulli bond percolation on $\bZ^d$. Therefore the time constant does not provide much information if $F(0)\geq p_c$.

In \cite[Eq.~3]{zhang1995supercritical}, Y. Zhang introduced the following random variable
\[
\rho(F) = \lim_{n\to\infty} T(\vzero,\partial B(n))\ ,
\]
where $B(n) = \{x\in\bZ^2:\|x\|_{\infty} \leq n\}$, $\partial B(n) = \{x\in\bZ^2:\|x\|_{\infty}=n\}$, and $\|\cdot\|_\infty$ is the sup-norm. By monotonicity, $\rho(F)$ exists almost surely. It was shown in \cite[p.~254]{zhang1995supercritical} that if $F(0)>p_c$ and $t_e$ has all moments, then for any $m\in\bN$, one has $\E\rho^m(F)<\infty$, and hence in particular $\rho(F)<\infty$ almost surely. Also, it is easy to see that if $F(0)<p_c$, then $\rho(F)=\infty$ almost surely. Then a natural question arises: how about $F(0)=p_c$?

In \cite{zhang1999double}, Zhang proved that for $d=2$, it is possible to have $\rho(F)<\infty$ or $\rho(F)=\infty$ almost surely when $F(0)=p_c$ (note that by the Kolmogorov zero-one law, either $\rho(F)<\infty$ almost surely or $\rho(F)=\infty$ almost surely). More specifically, he introduced the following two distributions. For $a>0$, set
\[
 F_a(x)= 
\begin{cases}
    1& \text{if } x^a > 1-p_c,\\
    x^a+p_c & \text{if } 0\leq x^a\leq 1-p_c,\\
    0 & \text{if } x<0,
\end{cases}
\]
and for $b>0$, set
\[
    G_b(x)= 
\begin{cases}
    1& \text{if } \exp(-1/x^b) > 1-p_c,\\
    \exp(-1/x^b)+p_c & \text{if } 0\leq \exp(-1/x^b)\leq 1-p_c,\\
    0 & \text{if } x<0.
\end{cases}
\]
Zhang showed in \cite[Theorem~8.1.1]{zhang1999double} that if $a$ is sufficiently small then $\rho(F_a)<\infty$ almost surely. He also made the following conjecture (see \cite[p.~146]{zhang1999double}): 
\begin{conj}[Zhang]\label{conj: Zhang}
The quantity $\sup\{a>0: \rho(F_a)<\infty\}$ is finite.
\end{conj}
\noindent
Moreover, Zhang showed in \cite[Theorem~8.1.3]{zhang1999double} that if $b>1$, then $\rho(G_b)=\infty$ almost surely.

 
The critical case of first-passage percolation is quite different from the standard one and requires different techniques. For example, the model is expected to retain rotational invariance in the limit \cite{yao2014law}, whereas the usual first-passage model has lattice dependent and distribution dependent asymptotics. For this reason, analysis of the critical case relies on detailed estimates from critical and near-critical percolation (for instance, see \cite{sapozhnikov2011incipient,jarai2003invasion,van2007size}). The main new insight of our work is that the behavior of passage times is closely related to a ``greedy'' growth algorithm called invasion percolation, and that optimal paths constrained to lie in the invasion cluster have the correct first-order growth. This relation allows us to derive necessary and sufficient conditions on the edge-weight distribution to have diverging mean or variance for passage times (Theorems~\ref{thm:mean} and \ref{thm:variance}), and these results can be seen as finer versions of Kesten's condition \eqref{eq: old_kesten} for $\mu=0$. Furthermore, we can derive a type of universality: for any edge weights for which the passage-time variance diverges, one has Gaussian fluctuations (see Theorem~\ref{thm:limit-theorem}).



Constants in this paper may depend on the distribution function $F$ and other fixed parameters such as $\eta$, $r$ and $\lambda$. However, constants do not depend on $k$ or $n$. We use $C_1,C_2, \ldots$ to denote temporary constants whose meaning may vary, while we use notation like $K_{\ref{lem:moments-t-k}}$ to denote the permanent constants. For example, $K_{\ref{lem:moments-t-k}}$ denotes the constant in Lemma \ref{lem:moments-t-k}.

\subsection{Main results}
In this paper, we will give an exact criterion for $\rho(F)<\infty$ (see Corollary~\ref{cor:rho-f} below) and consequently provide a negative answer to Conjecture~\ref{conj: Zhang}. Furthermore, we will derive limit theorems for the sequence $(T(\vzero ,\partial B(n)))_{n \geq 1}$. From now on, suppose that $d=2$ and that $F(0)=p_c$. Furthermore, define
\[
F^{-1}(t) = \inf\{x: F(x) \geq t\} \qquad \text{ for }  t > 0
\]
and
\begin{equation}
	\label{eqn:def-eta0}
	\eta_0 := \sup\set{\eta \geq 0 : \E [t_e^{\eta/4}] < \infty}.
\end{equation}

\subsubsection{Behavior of the mean}

We begin with bounds on $\E T(\vzero,\partial B(n))$.

\begin{thm}
	\label{thm:mean}
	(i) Assume that $\eta_0 > 1$. There exists $C_1 = C_1(F)> 0$ such that
	\begin{equation*}
		\E T(\vzero,\partial B(2^n)) \leq C_1 \sum_{k=2}^{n}F^{-1}(p_c+2^{-k}) \qquad \text{ for } n \geq 2.
	\end{equation*}	
	(ii) There exists $C_2 = C_2(F) > 0$ such that 
	\begin{equation*}
		\E T(\vzero,\partial B(2^n)) \geq C_2\sum_{k=2}^{n}F^{-1}(p_c+2^{-k}) \qquad \text{ for } n \geq 2.
	\end{equation*}
\end{thm}
%
\begin{rem}
Note that $\eta_0 > 1$ if and only if $\E Y^p<\infty$ for some $p>1$, where $Y$ is the minimum of four i.i.d. random variables distributed as $t_e$. The moment condition in Theorem~\ref{thm:mean} is nearly optimal since, if $\E Y = \infty$ then, by bounding $T(\vzero,\partial B(2^n))$ below by the minimum of the $4$ edge-weights on edges incident to $\vzero$, one has $\E T(\vzero,\partial B(2^n)) = \infty$ for $n\geq0$.
\end{rem}

\begin{rem}
\label{rem:point-to-point}
The above theorem concerns the passage time from the point $\vzero$ to the set $\partial B(n)$. In Section~\ref{sec: limit-0-to-x}, we derive asymptotics for point-to-point passage times $\E T(\vzero ,x)$ for $x \in \bZ^2$.
\end{rem}


As a corollary, we have an exact criterion for finiteness of $\rho(F)$.
\begin{cor}
\label{cor:rho-f}
For any $F$, one has $\rho(F)<\infty$ almost surely if and only if $\sum_{n=2}^\infty F^{-1}(p_c+2^{-n})<\infty$.
\end{cor}

We will now apply the above results to $F_a$ and $G_b$, the distributions defined by Zhang. The proof follows by a direct computation and the previous corollary.
\begin{cor}
	\label{cor:zhang}
The following statements hold.
\begin{enumerate}
\item $\rho(F_a)<\infty$ almost surely for any $a>0$, and so $\sup\{a>0: \rho(F_a)<\infty\} = \infty$. In particular, Conjecture~\ref{conj: Zhang} is false. 
\item $\rho(G_b)=\infty$ almost surely if and only if $b\geq 1$.
\end{enumerate}
\end{cor}

\begin{rem}
Zhang asked in \cite[p.~145]{zhang1999double} if, under the assumption $\E t_e^m<\infty$ for all $m\in\bN$, does $\rho(F)<\infty$ almost surely imply that $\E\rho(F)<\infty$? The answer is yes by combining all the above results.
\end{rem}


\subsubsection{Behavior of the variance and limit theorems}

Now we consider $\var(T(\vzero,\partial B(2^n)))$.
\begin{thm}\label{thm:variance} 
	Assume that $\eta_0>2$.\\
	(i) There exists $C_3 = C_3(F) > 0$ such that
	\begin{equation*}
		\var(T(\vzero,\partial B(2^n))) \leq C_3 \sum_{k=2}^{n}[F^{-1}(p_c+2^{-k})]^2 \qquad \text{ for } n \geq 2.
	\end{equation*}	
	(ii) There exists $C_4 = C_4(F) > 0$ such that 
	\begin{equation*}
		\var(T(\vzero,\partial B(2^n))) \geq C_4\sum_{k=2}^{n}[F^{-1}(p_c+2^{-k})]^2 \qquad \text{ for } n \geq 2.
	\end{equation*}
\end{thm}

By Corollary \ref{cor:rho-f}, when $\sum_{k=2}^\infty F^{-1}(p_c+2^{-k}) = \infty$ we have $T(\vzero, \partial B(n)) \convas \infty$ as $n \to \infty$, . The next theorem gives more information about the limit of $T(\vzero, \partial B(n))$ in this case.

\begin{thm}\label{thm:limit-theorem}
	Suppose $\sum_{k=2}^\infty F^{-1}(p_c+2^{-k}) = \infty$ and $\eta_0 > 2$.\\
	(i) If $\sum_{k=2}^\infty [F^{-1}(p_c+2^{-k})]^2 < \infty$, then there is a random variable $Z$ with $\E Z = 0$ and $\E Z^2 < \infty$ such that as $n \to \infty$
	\begin{equation*}
		T(\vzero,\partial B(n)) - \E T(\vzero,\partial B(n)) \to Z \quad \mbox{ a.s. and in $L^2$}.
	\end{equation*}
	(ii) If $\sum_{k=2}^\infty [F^{-1}(p_c+2^{-k})]^2 = \infty$, then as $n \to \infty$
	 \begin{equation*}
		\frac{T(\vzero,  \partial B(n)) - \E T(\vzero,  \partial B(n))}{{[\var(T(\vzero,  \partial B(n)))]}^{1/2}}
		 \convd N(0,1).
	\end{equation*}
\end{thm}

\begin{rem}
As in the case of Theorem~\ref{thm:mean}, in Section~\ref{sec: limit-0-to-x}, we derive versions of the variance asymptotics and limit theorems for point-to-point passage times $T(\vzero,x)$ for $x \in \bZ^2$. See Corollaries~\ref{cor:variance-p2p} and \ref{cor:limit-theorem-p2p}.
\end{rem}

\subsection{Relations to previous work}

First-passage percolation has been studied since its introduction by Hammersley and Welsh \cite{hammersley1965first} in the '60s, but most work has focused on the non-critical case, where $F(0)<p_c$. There, the passage time from $\vzero$ to a vertex $x$ grows linearly in $x$, and many results have been proved, including shape theorems, large deviations, concentration inequalities and moment bounds. We refer the reader to the surveys \cite{grimmett2012percolation, blair2010first}. The supercritical case, where $F(0)>p_c$ is easier to analyze, since there is almost surely an infinite cluster of edges with passage time 0, and so distant vertices need only to travel to the infinite cluster to reach one-another. This produces passage times $T(\vzero,x)$ that are of order one as $x \to \infty$.

The critical case, where $F(0)=p_c$, is considerably more subtle. It is expected (though only proved in two dimensions or high dimensions) that there is no infinite cluster of $p_c$-open edges (that is, edges with passage time 0 in this case). However, clusters of $p_c$-open edges occur on all scales, giving, for example, infinite mean size for the $p_c$-open cluster of the origin. This means that two distant points can be connected by a path which uses mostly zero-weight edges, and this path may be able to find lower and lower edge weights as it moves further into the bulk of the system. Therefore to characterize passage times, one should understand the balance between the number of edges on each scale with low weights and the number of paths that can access them.

Kesten proved in \cite[Theorem~6.1]{kesten1986aspects} that the time constant $\mu$ is zero in the critical case, implying that $T(\vzero, x) = o(\|x\|)$ as $x \to \infty$. This result was sharpened by L. Chayes \cite[Theorem~B]{chayes1991critical}, who showed that for any $\delta>0$, $\lim_{n \to \infty} T(\vzero, n\mathbf{e}_1) / n^\delta = 0$ almost surely. In \cite[Remark~3]{kesten1993speed}, Kesten claimed that in fact Chayes's argument can be extended to $T(\vzero, n\mathbf{e}_1) \leq \exp(C\sqrt{\log n})$ for large $n$ almost surely. These results go some way to quantify asymptotics of the passage time in the critical case for general dimension.

More progress has been made in the critical case in two dimensions, due to a more developed theory of Bernoulli percolation on planar lattices. It was shown by Chayes-Chayes-Durrett in \cite[Theorem~3.3]{chayes1986critical} that if $t_e$ is Bernoulli (0 or 1 with probability 1/2) then the expected passage time grows logarithmically, obeying $\mathbb{E}T(\vzero, n\mathbf{e}_1) \asymp \log n$. In this Bernoulli case, the passage time between $\vzero$ and $x$ can be represented as the maximum number of disjoint $p_c$-closed circuits separating $\vzero$ and $n\mathbf{e}_1$, as every $p_c$-closed edge on a geodesic contributes passage time 1. Recently, Yao \cite{yao2014law} has shown a law of large numbers on the triangular lattice, using the conformal loop ensemble of Camia and Newman.

Our work was motivated by that of Zhang in '99, who showed that critical FPP can display ``double behavior." That is, he showed that there exist distributions $F$ with $F(0)=p_c$ for which the passage time $T(\vzero, \partial B(n))$ diverges as $n\to\infty$, and those for which the passage time remains bounded. Intuitively, bounded passage times come from those distributions which have significant mass near zero, so that long paths can find more and more low weights as they move away from $\vzero$, producing infinite paths with finite passage time. Zhang asked many questions about this case, in particular which distributions have which of the two behaviors. One main point of our work is Theorem~\ref{cor:rho-f}, which gives an exact criterion that this passage time remains bounded if and only if $\sum_k F^{-1}(p_c + 2^{-k}) < \infty$. Our proof involves a new relation to a model called invasion percolation, and it turns out that optimal paths in the invasion cluster have passage time of the same order as geodesics in FPP. (See the next section for more details.) This theorem allows us to answer Zhang's questions in the two-dimensional case.

The other motivation for our work is that of Kesten and Zhang in '97. They also considered the critical case in two dimensions and proved central limit theorems for $T(\vzero, \partial B(n))$ for a certain class of distributions. Precisely, they showed that if $\mathbb{E}t_e^\delta<\infty$ for some $\delta>4$, $F(0)=p_c$, and there exists a constant $C_0>0$ such that $F(C_0) = p_c$, then  the sequence $T(\vzero, \partial B(n))$ satisfies a Gaussian central limit theorem: there exists a sequence $\gamma_n$ such that
\[
C_1(\log n)^{1/2} \leq \gamma_n \leq C_2(\log n)^{1/2}
\]
and
\[
\frac{T(\vzero, \partial B(n)) - \mathbb{E}T(\vzero, \partial B(n))}{\gamma_n} \Rightarrow N(0,1).
\]
It is important to notice that the condition $F(C_0)=p_c$ gives a positive lower bound for the passage time of non-zero weight edges. Kesten and Zhang do not address any distributions with mass near zero, though they do remark about the double behavior of such distributions.

 The second part of our paper, on limit theorems and variance estimates, completes the picture started by Kesten and Zhang. Theorems~\ref{thm:variance} and \ref{thm:limit-theorem} (ii) require only that $\sum_k (F^{-1}(p_c+2^{-k}))^2 = \infty$ and a weak moment condition on $t_e$ (lower than that of Kesten and Zhang) to deduce that the variance of $T(\vzero, \partial B(n))$ diverges and that a Gaussian CLT holds. This result on the CLT shows that in the critical case, no other limiting behavior is possible, in contrast to the subcritical case, where the variance is expected to be of order $n^{2/3}$ with a non-Gaussian limiting distribution (see \cite{kardar1986dynamic}). Theorem~\ref{thm:limit-theorem} (i) also addresses the intermediate case, where the mean of $T(\vzero, \partial B(n))$ diverges but the variance converges. Here, the centered sequence is tight and converges to a non-trivial limit. We do not know an exact form for this limit, and it is unlikely to be explicit since its variance depends heavily on weights of edges near the origin.


\section{Setup for the proof}
\label{sec:setup-for-proof}

Zhang's proof in \cite[Theorem~8.1.1]{zhang1999double} that $\rho(F_a)$ has all moments used a comparison to a near-critical percolation model introduced in \cite{chayes1987inhomogeneous} by Chayes-Chayes-Durrett. Their model is a version of an incipient infinite cluster, a term used by physicists to describe large (system-spanning) percolation clusters at criticality. We will, however, need finer asymptotics that are obtained by comparison with a different near-critical model, invasion percolation. Though it has no parameter, it tends on large scales to resemble Bernoulli percolation at criticality. We describe the model of invasion percolation in Section \ref{sec:invasion-percolation}. We also recall some known facts about Bernoulli percolation in Section \ref{sec:correlation-length}.

We will couple the first-passage percolation model on $(\bZ^2,\cE^2)$ with invasion percolation and Bernoulli percolation. To describe the coupling, we consider the probability space $(\Omega, \cF, \pr)$, where $\Omega = [0,1]^{\cE^2}$, $\cF$ is the cylinder sigma-field and $\pr = \prod_{e \in \cE^2} \mu_e$, where each $\mu_e$ is an uniform distribution on $[0,1]$. Write $\omega=(\omega_e)_{e \in \cE^2} \in \Omega$.  Define the edge weights as $t_e = F^{-1}(\omega_e)$ for $e \in \cE^2$.

\subsection{Invasion percolation}
\label{sec:invasion-percolation}


If an edge $e$ has endpoints $e_x$ and $e_y$, we write $e=\{ e_x,e_y\}$. For an arbitrary subgraph $G=(V,E)$ of $(\bZ^2,\mathcal{E}^2)$, define the edge boundary $\Delta G$ by
\[
\Delta G = \{e \in\mathcal{E}^2: e\not\in E, e_x\in V\text{ or } e_y\in V\}.
\]
Define a sequence of subgraphs $(G_n)_{n=0}^\infty$ as follows. Let $G_0 = (\{0\},\emptyset)$. If $G_i=(V_i,E_i)$ is defined, we let $E_{i+1} = E_i \cup \{e_{i+1}\}$, where $e_{i+1}$ is the edge with $\omega_{e_{i+1}}= \min\{\omega_e: e\in\Delta G_i\}$, and let $G_{i+1}$ be the graph induced by $E_{i+1}$. The graph $I:=\bigcup_{i=0}^\infty G_i$ is called the \emph{invasion percolation cluster} (at time infinity).


Invasion percolation is coupled with the first-passage percolation model since we have defined $t_e = F^{-1}(\omega_e)$. They can also be coupled with Bernoulli percolation as follows. For each $e \in \mathcal{E}^2$ and $p \in [0,1]$, we say that $e$ is $p$-open in $\omega$ if $\omega_e \leq p$ and otherwise we say that $e$ is $p$-closed. If there is a $p$-open path from a vertex set $A$ to a vertex set $B$ then we write that $A \leftrightarrow B$ by a $p$-open path. The collection of $p$-open edges has the same distribution as the set of open edges in Bernoulli percolation with parameter $p$. 

To use this coupling, we need the notion of the dual graph. Let $(\bZ^2)^* = (1/2,1/2)+\bZ^2$ and $(\mathcal{E}^2)^* = (1/2,1/2)+\mathcal{E}^2$. For $x\in\bZ^2$, we write $x^* = (1/2,1/2)+x$. For $e\in\mathcal{E}^2$, we denote the its endpoints (left respectively right or bottom respectively top) by $e_x$, $e_y\in\bZ^2$. The edge $e^*=\{ e_x+(1/2,1/2),e_y-(1/2,1/2) \}$ is called the dual edge to e and its endpoints (bottom respectively top or left respectively right) are denoted by $e_x^*$ and $e_y^*$. For $A\subset\bZ^2$, $A^*$ is defined to be $(1/2,1/2)+A$. An edge $e^*$ is declared to be $p$-open in $\omega$ when $e$ is, and $p$-closed otherwise.

We note the following relations between invasion percolation and Bernoulli percolation:
\begin{itemize}
\item With probability one, if $x$ is a vertex of $I$ and $y \leftrightarrow x$ by a $p_c$-open path, then $y \in I$.

\textbf{Proof: }
If $y$ is not in $I$ then we can find $e \in \Delta I$ (on a $p_c$-open path from $x$ to $y$) such that $e$ is $p_c$-open. But then $e \in \Delta G_n$ for all large $n$. By the definition of the invasion algorithm, this means that for large $n$, each edge added to the invasion is $p_c$-open, and from this we can build an infinite $p_c$-open path. This contradicts the fact that there is almost surely no infinite $p_c$-open cluster \cite[Theorem~1]{kesten1980critical}.
\qed\\
\item For $n \geq 0$, let $\hat p_n$ be defined as
\begin{equation}\label{eqn:def-hat-pn}
\hat p_n = \sup\{\omega_e : e \in I \cap E(B(2^n))^c\}\ ,
\end{equation}
where $E(V)$ is the set of edges with both endpoints in $V$. Then
\begin{equation}\label{eqn:old-implication}
\hat p_n > p \Rightarrow A_{n,p} \text{ occurs}\ ,
\end{equation}
where
\[
A_{n,p} = \{\exists ~p\text{-closed dual circuit around the origin with diameter at least }2^n\}\ .
\]
Here the diameter of a set $X$ is $\sup\{\|x-y\|_\infty : x,y \in X\}$.

\textbf{Proof: }
Take $e \in I\cap E(B(2^n))^c$ with $\omega_e>p$. At the moment $k$ that $e$ is added to the invasion cluster, the graph $G_k$ has edge boundary all of whose edges are $\omega_e$-closed, and so are $p$-closed. However, from the edge boundary, we can extract a dual circuit around 0 that contains $e^*$, by \cite[Proposition 11.2]{grimmett1999percolation}. This circuit then has diameter at least $2^n$.
\qed\\
\end{itemize}

\subsection{Correlation length}
\label{sec:correlation-length}
A central tool used to study invasion percolation is correlation length, and we take the definition from \cite[Eq.~1.21]{kesten1987scaling}. For $m,n\in\bN$ and $p\in (p_c,1]$, let
\[
\sigma(n,m,p) = \pr(\text{there is a } p\text{-open left-right crossing of }[0,n]\times[0,m]),
\] 
where a $p$-open left-right crossing of $[0,n]\times[0,m]$ means a path $\gamma$ in $[0,n]\times[0,m]$ with all edges $p$-open which joins some vertex on $\{0\}\times[0,m]$ to some vertex on $\{n\}\times[0,m]$.
For $\epsilon>0$ and $p> p_c$, we define
\[
L(p,\epsilon) = \min\{n \geq 1: \sigma(n,n,p) \geq 1-\epsilon\}.
\]
$L(p,\epsilon)$ is called the correlation length. It is known (see \cite[Eq.~1.24]{kesten1987scaling}) that there exists $\epsilon_1>0$ such that for all $0<\epsilon,\epsilon'\leq \epsilon_1$, the ratio $L(p,\epsilon)/L(p,\epsilon')$ is bounded away from 0 and $\infty$ as $p\downarrow p_c$. We will write $L(p)=L(p,\epsilon_1)$ for simplicity. For $n\geq 1$, define 
\begin{equation}
	\label{eqn:def-pn}
	p_n = \min\{p: L(p)\leq n\}.
\end{equation}
 
We now note the following facts.

\begin{itemize}
\item By \cite[Eq. (2.10)]{jarai2003invasion} there exists $K_{\ref{eqn:jarai-correlation-length}} \in (0,1)$ such that for all $n \geq 1$ we have
 	\begin{equation} \label{eqn:jarai-correlation-length}
 		K_{\ref{eqn:jarai-correlation-length}}  n \leq L(p_n) \leq n.
 	\end{equation}
\item There exist $C_1, C_2 > 0$ such that for all $m,n\geq 1$,
\[
C_1\left|\log\frac{m}{n}\right|\leq\left|\log\frac{p_m-p_c}{p_n-p_c}\right|\leq C_2\left|\log\frac{m}{n}\right|.
\]
This is a consequence of \cite[Prop.~34]{nolin2008near} and a priori estimates on the four-arm exponent. In particular, putting $m=1$, there exist $ \delta_0 > \epsilon_0 > 0$ such that for $n \geq 2$
\begin{equation}
 	\label{eqn:466}
 	 \frac{1}{n^{\delta_0} } < p_n - p_c < \frac{1}{n^{\epsilon_0} }.
\end{equation}
We may and will always assume $\delta_0 >1$. 

\item From \cite[Eq.~9]{van2008random} and \eqref{eqn:old-implication}, There exist $K_{\ref{eqn:outlet}.1}, K_{\ref{eqn:outlet}.2} > 0$ such that for all $p > p_c$ and $n \geq 1$,
\begin{equation}
\label{eqn:outlet}
 \pr(\hat p_n > p) \leq \pr(A_{n,p}) \leq K_{\ref{eqn:outlet}.1} \exp\left(-\frac{K_{\ref{eqn:outlet}.2} 2^n}{L(p)} \right).
\end{equation}
\item By the RSW theorem (see \cite[Section~11.7]{grimmett1999percolation}), there exists $K_{\ref{eqn:5-26}}>0$ such that for all $k\in\bN$,
\begin{equation}
\label{eqn:5-26}
\pr(\text{there exists a }p_{2^k}\text{-closed dual circuit around } \vzero \text{ in }B(2^k)^* \setminus B(2^{k-1})^* )\geq K_{\ref{eqn:5-26}}.
\end{equation}
\end{itemize}

\subsection{Sketch of proofs} 

The central tool used to prove our theorems is Lemma~\ref{lem:moments-t-k}, which is a moment bound on annulus passage times. We first describe the idea of its proof. Consider all paths between $\vzero$ and $\partial B(2^{n+1})$ which lie in the invasion cluster $I$ and $B(2^{n+1})$. Let $\gamma_n$ be such a path which minimizes the passage time. Then Lemma~\ref{lem:moments-t-k} gives an upper bound on the $r$-th moment of the sum of edge-weights for edges in $\gamma_n$ which lie in any annulus $B(2^{k+1}) \setminus B(2^k)$ (that is, $\mathbb{E}T_k^r(\gamma_n)$, where $T_k(\gamma_n)$ is defined in \eqref{eqn:def-T-k}).

The passage time of $\gamma_n$ gives an upper bound for the passage time from $\vzero$ to $\partial B(2^n)$, and $\gamma_n$ is in many ways is a nicer path than the actual geodesic using the edge-weights $(t_e)$. Once the invasion has reached the boundary of $B(2^k)$, all of its edges (from that point on) are likely to be nearly $p_{2^k}$-open (that is, $\hat p_k$ from \eqref{eqn:def-hat-pn} is of order $p_{2^k}$), and so the edges in $\gamma_n$ outside of $B(2^k)$ will have passage time bounded above by $F^{-1}(p_{2^k})$. By bounding $p_{2^k}$ above with \eqref{eqn:466}, each edge has passage time bounded by $a_k$, where $a_k$ is defined in \eqref{eq: a_k_def}. Unfortunately we only know this behavior of $\hat p_k$ with high probability, so we need to decompose the probability space over different values of $\hat p_k$ using an idea of A. J\'arai \cite{jarai2003invasion}.

The above heuristic gives
\[
T_k(\gamma_n) \lesssim a_k \#\{e \in \gamma_n \cap (B(2^{k+1}) \setminus B(2^k)) : e \text{ is } p_c\text{-closed}\}.
\]
The reason is that the only edges which contribute to the passage time of $\gamma_n$ are those which are $p_c$-closed. In Lemma~\ref{lem:main-lemma}, we show that each such edge has ``4-arms.'' That is, they have the property that (a) their weight is between $p_c$ and $p_{2^k}$, (b) they have two disjoint $p_{2^k}$-open arms to distance $2^{k-1}$ and (c) they have two disjoint $p_c$-closed arms to distance $2^{k-1}$. Fortunately all moments of the number of such points in an annulus were bounded in the study of invasion percolation in \cite{damron2011outlets} (see Lemma~\ref{lem:498} below), so we can conclude.

\subsubsection{Idea of the proof of Theorem~\ref{thm:mean}}
The proof of (ii) follows that of Zhang \cite[Theorem~8.1.2]{zhang1999double}. The proof of (i) follows immediately from Lemma~\ref{lem:moments-t-k}. Indeed, to find the upper bound for $\E T(\vzero, \partial B(2^n))$, we simply use the inequality
\[
T(\vzero, \partial B(2^n)) \leq T(\gamma_n) = \sum_{k=-1}^n T_k(\gamma_n),
\]
where, as above, each $T_k(\gamma_n)$ is the time that $\gamma_n$ spent in the annulus $B(2^{k+1})\setminus B(2^{k})$. Applying the annulus moment bounds from Lemma~\ref{lem:moments-t-k} gives (i).

\subsubsection{Idea of the proof of Theorem~\ref{thm:variance} and \ref{thm:limit-theorem}}
To study the variance and limit theorems, we follow the strategy of Kesten-Zhang \cite{kesten1997central}. Instead of dealing with $\var(T(\vzero,\partial B(2^n)))$ directly, we consider $\var(T(\vzero,\mathcal{C}_n))$, where $\mathcal{C}_n$ is the innermost $p_c$-open circuit in an annulus $B(2^{m+1}) \setminus B(2^m)$ for $m \geq n$ surrounding $\vzero$. It can be shown that these two variances are closed to each other. The variance bounds for $T(\vzero, \mathcal{C}_n)$ are stated in Theorem~\ref{thm:variance-zero-circuit} and the CLT is stated in Theorem~\ref{thm:clt-zero-circuit}.

If we write $T(\vzero,\mathcal{C}_n) - \E T(\vzero, \mathcal{C}_n)$ as a sum of martingale differences 
\[
\Delta_k = \mathbb{E}[T(\vzero, \mathcal{C}_n) \mid \mathcal{F}_k] - \mathbb{E}[T(\vzero, \mathcal{C}_n) \mid \mathcal{F}_{k-1}]
\]
over a filtration $(\mathcal{F}_k)$, then $\var (T(\vzero,\mathcal{C}_n)) = \sum_{k=0}^n \E \Delta_k^2$, so it suffices to bound the $\E \Delta_k^2$'s. The idea of Kesten-Zhang was to take $\mathcal{F}_k$ to be generated by the edge-weights for edges on and in the interior of $\mathcal{C}_k$, and they proved an alternate representation for such $\Delta_k$'s (see Lemma~\ref{lem:delta-p-decomposition} (ii)). With this same choice, we can use the moment bounds in Lemma~\ref{lem:moments-t-k} to prove moment bounds on the $\Delta_k$'s in Lemma~\ref{lem:delta-p-upper-bound}. We emphasize that in fact, as a result of the representation in Lemma~\ref{lem:delta-p-decomposition} (ii), $\Delta_k$ does not depend on $n$. 

To prove the CLT for $T(\vzero, \partial B(n))$, we apply McLeish's CLT, which we state as Theorem~\ref{thm:mart-clt-mcleish}. We verify its conditions using the fact that the $\Delta_k$'s can be shown to be strongly mixing (see Lemma~\ref{lem:strong-mixing-coefficient}),  Along with moment bounds for the $\Delta_k$'s given in Lemma~\ref{lem:delta-p-upper-bound} (which use our main moment bounds on annulus times), we can conclude the CLT in the case that the variance of $T(\vzero, \partial B(n))$ diverges. This will prove (ii) in Theorem~\ref{thm:limit-theorem}. For (i), if the variance does not diverge, then by the martingale convergence theorem, $T(\vzero,\mathcal{C}_n) - \E T(\vzero,\mathcal{C}_n)$ will converge to some random variable $Z$. Using a stronger comparison to $T(\vzero, \partial B(n))$ given in Lemma~\ref{lem:error-bn-cq} allows to complete the proof.

\section{Moment bounds for annulus times}
\label{sec:moment-bounds}
In this section, we prove the main lemma of the paper, Lemma~\ref{lem:moments-t-k}. It serves to bound certain annulus passage times $T_k(\gamma_n)$ through the invasion cluster.

Recall that we denote $I$ as the invasion percolation cluster using the weights $\set{\omega_e}_{e \in \cE^2}$. Define $\cE_{-1}:= E(B(1))$ and  $\cE_n:= E(B(2^{n+1})) \setminus E(B(2^{n}))$ for $n \geq 0$. Note that $|\cE_{-1}| = 12$ and $|\cE_{n}| = 24 \cdot 4^n + 4 \cdot 2^n$ for $n \geq 0$.  For any path $\gamma$, define for $ k \geq -1 $
\begin{equation} \label{eqn:def-T-k}
    T_k(\gamma) := \sum_{ e \in \gamma \cap \cE_k} t_e.
\end{equation}

 For $n \geq -1$, let $\gamma_n$ be a path such that 
 \begin{equation*}
 	T(\gamma_n) = \inf \set{ T(\gamma) : \gamma \text{ is a path between $\vzero$ and $\partial B(2^{n+1})$ and $\gamma \subset B(2^{n+1})\cap I$}  }.
 \end{equation*}
Note that $T(\gamma_n) = \sum_{k=-1}^{n} T_k(\gamma_n)$.
Recall $\epsilon_0$ from \eqref{eqn:466}. For simplicity of notation, define
 \begin{equation}
 \label{eq: a_k_def}
 	a_k := F^{-1}(p_c + 2^{-\epsilon_0  k/2}), \mbox{ for } k \in \bN.
 \end{equation}
 Note that $a_k$ is only defined when the argument of $F^{-1}$ is strictly less than $1$, and this will be guaranteed by the condition $k \geq k_0$ in the lemma below.

 The main goal of this section is to prove the following lemma.
 \begin{lem}
 	\label{lem:moments-t-k} 
	Recall the definition of $\eta_0$ from \eqref{eqn:def-eta0} and suppose $\eta_0 > 1$.
	
 	(i) For all $r \in [1,\eta_0)$ and integers $k \geq -1$, we have  $ \sup_{n \geq k} \E[ T_k^r(\gamma_{n}) ] < \infty$.
	

 	(ii) Given any $r \in [1,\infty)$ and $\lambda \in (0,\infty)$, there exist $k_0 = k_0(r,\lambda,F) > 0$ and $K_{\ref{lem:moments-t-k}} = K_{\ref{lem:moments-t-k}}(r, \lambda, F) > 0$, such that for all $n - 1 \geq k  \geq k_0$ we have
 	\begin{equation*}
 		\E[ T_k^r(\gamma_{n}) ] \leq K_{\ref{lem:moments-t-k}}(a_k^r + e^{-\lambda k}).
 	\end{equation*}
 \end{lem}
 
 \begin{rem}
 	To prove Theorem \ref{thm:mean}, it is sufficient to use the above lemma with $r = 1$. Here we prove it in the general form for future use in Section \ref{sec:proof-variance}.
 \end{rem}
 
We begin with a definition from \cite{damron2011outlets}. For $m_1,m_2\geq 1$, $p\in(p_c,1]$, and $e \in \mathcal{E}^2$, let $A_e(m_1,p)$ be the event that
\begin{enumerate}
\item[(a)] $e$ is connected to $\partial B(e_x,m_1)$ by two vertex disjoint $p$-open paths, 
\item[(b)] $e^*$ is connected to $\partial B(e_x,m_1)^*$ by two vertex disjoint $p_c$-closed dual paths, and
\item[(c)] $\omega_e \in (p_c,p]$.
\end{enumerate}
Here $\partial B(e_x,m_1) = e_x + \partial B(m_1)$. Let $N(m_1,m_2,p)$ be the number of edges $e$ in $E(B(2m_2)) \setminus E(B(m_2))$ such that $A_e(m_1,p)$ occurs; that is,
\[
N(m_1,m_2,p) = \sum_{e \in E(B(2m_2)) \setminus E(B(m_2))} \ind_{A_e(m_1,p)}.
\]

\begin{lem}
 	\label{lem:main-lemma}
	Let $\hat p_k$ be as in \eqref{eqn:def-hat-pn}. For all $p > p_c$ and $1 \leq k \leq n-1$ we have
 	\begin{equation*}
 		T_k(\gamma_{n})\ind\set{\hat p_k \leq p} \leq N(2^{k-1}, 2^k, p) \cdot F^{-1}(p).
 	\end{equation*}
 \end{lem}
\textbf{Proof: }
Suppose $\hat p_k \leq p$ for some $p > p_c$. Define, for $n \geq 1$ and $1 \leq k \leq n-1$
\begin{align*}
T_{k,n}' = \# \{ e \in \gamma_n\cap \mathcal{E}_k : \omega_e > p_c\}\ .
\end{align*}
Since $\hat p_k \leq p$ and $\gamma_n \subset I$, we have $T_k(\gamma_n) \leq T_{k,n}' F^{-1}(p)$. Then it is sufficient to show
\begin{equation}\label{eq: claim_1}
T_{k,n}' \leq N(2^{k-1},2^k,p).
\end{equation}

Let $e\in\gamma_n\cap \mathcal{E}_k$ be $p_c$-closed. As $\gamma_n\subset I$ and $\hat p_k \leq p$, $e$ is $p$-open. Note that there exist disjoint paths $\gamma_{n,1}$, $\gamma_{n,2}\subset\gamma_n$ such that $\gamma_{n,1}$ is a $p$-open path joining $e_x$ to $\partial B(e_x,2^{k-1})$ and $\gamma_{n,2}$ is a $p$-open path joining $e_y$ to $\partial B(e_x,2^{k-1})$. (This holds because $e_x$ is invaded but $0 \notin B(e_x,2^{k-1})$.) 


\begin{figure}[t]
\begin{center}
\vspace{1.3cm}

\includegraphics[width=3.5in]{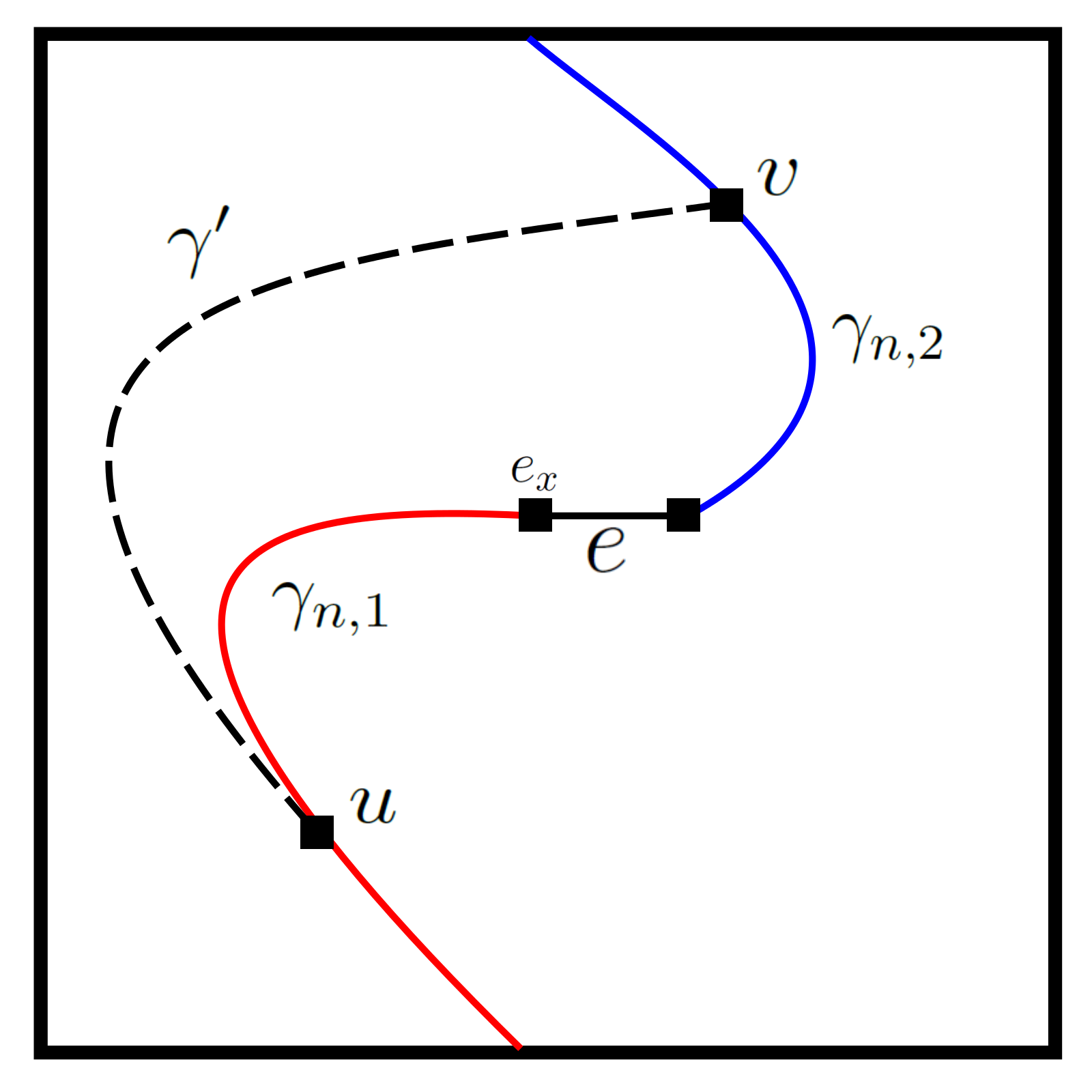}
\vspace{0.7cm}
\caption{Depiction of the proof of Lemma~\ref{lem:main-lemma}. The box shown is $B(e_x,2^{k-1})$. The path $\gamma'$ is $p_c$-open and connects vertices $u$ and $v$ on $\gamma_n$, but bypasses the edge $e$.}
\label{fig: fig_1}
\end{center}
\end{figure}


For an illustration of the following argument, see Figure~\ref{fig: fig_1}. If $\gamma_{n,1}\leftrightarrow\gamma_{n,2}$ by a $p_c$-open path $\gamma'$ in $B(e_x,2^{k-1})$, and if we let $u\in \gamma_{n,1}$ and $v\in\gamma_{n,2}$ be such that $u\leftrightarrow v$ via $\gamma'$, then every vertex in the path $\gamma'$ is in $I$ (see the first bulleted fact in Section~\ref{sec:invasion-percolation}). Therefore $\gamma'\subset I$. Now let $\gamma_n'$ be the path which connects $0$ and $u$ via $\gamma_n$, $u$ to $v$ via $\gamma'$ and $v$ to $\partial B(2^{n+1})$ via $\gamma_n$. Then $\gamma_n'$ is in $I$ and has at least one $p_c$-closed edge less (namely $e$) than $\gamma_n$. Furthermore, each $p_c$-closed edge of $\gamma_n'$ is a $p_c$-closed edge of $\gamma_n$, and this implies $T(\gamma_n') < T(\gamma_n)$, contradicting the minimality of $\gamma_n$. Hence $\gamma_{n,1}\not\leftrightarrow\gamma_{n,2}$ by a $p_c$-open path in $B(e_x,2^{k-1})$. Note that by duality, exactly one of the following will happen:

\begin{enumerate}
\item $e_x^*$ and $e_y^*$ are connected to $\partial B(e_x,2^{k-1})^*$ by two disjoint $p_c$-closed dual paths, which are also disjoint from $\gamma_{n,1}\cup \gamma_{n,2} \cup \set{e}$;
\item there is a $p_c$-open path connecting $\gamma_{n,1}$ and $\gamma_{n,2}$ in $B(e_x,2^{k-1})$.
\end{enumerate}

So the first event must happen and thus $A_e(2^{k-1},p)$ occurs. This completes the proof of Lemma \ref{lem:main-lemma}. \qed\\

Next we bound the moments of $N(2^{k-1}, 2^k, p)$ using a result from \cite[Lemma 5.1]{damron2011outlets}.
\begin{lem}
 	\label{lem:498}
 	 There exists $K_{\ref{lem:498}} > 0$ such that for all $p > p_c$, $L(p)< m_1 \leq m_2$ and integers $t \geq 1$
 	\begin{equation*}
 		\E[N^t(m_1,m_2,p)] \leq \E[N^t(L(p),m_2,p)] \leq t! \left( \frac{K_{\ref{lem:498}} m_2}{L(p)}\right)^{2t}.
 	\end{equation*}
\end{lem}
\textbf{Proof: } The first inequality immediately follows from the definition of $N(m_1,m_2,p)$. In \cite[Lemma 5.1]{damron2011outlets}, it was shown that there exists $C_1 > 0$ such that if $p>p_c$, $m'\leq L(p)$ and $m'\leq m_2$, then for all integers $t \geq 0$,
\begin{equation*}
\E [N^t(m',m_2,p)] \leq t!\left(C_1 \frac{m_2}{m'}\right)^{2t}.
\end{equation*}
Taking $m'=L(p)$ completes the proof.  \qed\\

The next lemma will be used to control moments of $T_k(\gamma_n)$ when $\hat p_k$ is large. Define

\begin{equation}
	\label{eqn:def-hat-t-k}
	\hat t_k := F^{-1}(\hat p_k).
\end{equation}

\begin{lem}
	\label{lem:hat-t-k-moments}
	Suppose $\E [t_e^\eta] < \infty$ for some $\eta > 0$. Define $c_{-1} = 4$ and $c_k := 2^{k+1} +4$ for $k \geq 0$. Then for all integers $k \geq -1 $ and $r \in (0, c_k \eta)$, one has $\E [\hat t_k^r] < \infty$. In particular, for any fixed $r > 0$, there exists $K_{\ref{lem:hat-t-k-moments}} = K_{\ref{lem:hat-t-k-moments}}(r,\eta,F)$ such that for all integers $k > \log({r}/{\eta})/ \log 2$, we have
\begin{equation*}
	\E [\hat t_k^r] \leq K_{\ref{lem:hat-t-k-moments}}.
\end{equation*}	
\end{lem}
\textbf{Proof: }  Note that $t \geq F^{-1}(F(t))$ for all $t \geq 0$. Then we have
\begin{equation}
	\label{eqn:974}
	\pr(\hat t_k > t) \leq \pr ( \hat t_k > F^{-1}(F(t))) \leq \pr( \hat p_k \geq F(t)).
\end{equation}
In order to bound the tail probability of $\hat t_k$, we need to bound $\pr(\hat p_k > p)$ when $p$ is close to one. By \eqref{eqn:old-implication}, for any $k \geq -1$, $\hat p_k > p$ implies that there exists a $p$-closed \emph{dual} circuit surrounding the origin with diameter at least $\lfloor 2^k+1 \rfloor$. Such a dual circuit must have length at least $2 \lfloor 2^k+1 \rfloor + 2  = c_k$, for $k \geq -1$.  For any even $m \geq 4$, observe that since dual circuits around the origin with length $m$ must intersect the line $\set{(x,0): x \in (-1, m/2-1)}$, the total number of such circuits is bounded by $\frac{m}{2} \cdot 3^m$. Each of these dual circuit is $p$-closed with probability $(1-p)^m$. Therefore when $p \in [5/6,1)$ we have
\begin{align*}
	 \pr( \hat p_k \geq p) \leq \sum_{m = c_k}^\infty \frac{m 3^m}{2} \cdot (1-p)^m \leq \sum_{m = c_k}^\infty \frac{m}{2^{(1-\alpha)m}}  (3(1-p))^{\alpha m} 
\end{align*}
where the second inequality uses the facts that $3(1-p) \leq 1/2$ and the value of $\alpha \in (0,1)$ will be specified later. Define $C_1 = C_1(\alpha):= \max_{m \geq 4} \set{m 2^{-(1-\alpha)m}}/(1 - 2^{-\alpha})$ and $C_2:= (3 \E t_e^{\eta})^{1/\eta}$. Combining \eqref{eqn:974} and the above bound, when $t \geq C_3 := F^{-1}(5/6)/C_2$ we have $F(C_2 t) \geq 5/6$ and
\begin{equation*}
	\pr( \hat t_k > C_2 t ) \leq \pr ( \hat p_k \geq F(C_2 t)) \leq C_1(3 \pr(t_e > C_2 t))^{ c_k\alpha} \leq C_1 \bpare{ \frac{3 \E[t_e^\eta]}{(C_2 t)^{\eta}}}^{c_k \alpha } = \frac{C_1}{t^{ c_k \alpha \eta}}.
\end{equation*}
Since $r < c_k \eta$, taking $\alpha = \alpha_k := \frac{c_k \eta + r}{2c_k\eta} $ we have
\begin{align}
	\E\bbrac{\bpare{\frac{\hat t_k}{C_2}}^r} 
	=& \int_0^\infty r t^{r-1} \pr( \hat t_k \geq C_2 t ) dt  \nonumber\\
	\leq& \int_0^{1 \vee C_3} r t^{r-1} dt   + \int_{1 \vee C_3}^\infty r t^{r-1} \cdot \frac{C_1(\alpha_k)}{t^{c_k\alpha_k \eta}} dt \nonumber\\
	\leq&  (1 \vee C_3)^r + \frac{C_1(\alpha_k) r}{c_k \alpha_k \eta -r}. \nonumber
\end{align}
Therefore, using the relation $c_k \alpha_k \eta -r = (c_k \eta -r)/2$, we have
\begin{equation*}
	\E[\hat t_k^r] \leq (C_2 \vee F^{-1}(5/6))^r + \frac{2r C_1(\alpha_k) C_2^r}{c_k \eta -r}.
\end{equation*}
This proves the statement that $\E[\hat t_k^r] < \infty$.

Next, when $r >0$ and $k > \log({r}/{\eta})/ \log 2$, taking $\alpha := 1/2$ in the above proof, we have $c_k \alpha \eta - r \geq 2^{k} \eta - r  \geq 2^{k_1} \eta - r > 0$ where $k_1 := \lfloor\log({r}/{\eta})/ \log 2 \rfloor + 1$. Then
\begin{equation*}
	\E[\hat t_k^r] \leq C_2^r + (F^{-1}(5/6))^r + \frac{2r C_1(1/2) C_2^r}{2^{k_1} \eta -r},
\end{equation*}
which gives the expression of $K_{\ref{lem:hat-t-k-moments}}$. \qed\\

Now we are ready to prove Lemma \ref{lem:moments-t-k}.

 \textbf{Proof of Lemma \ref{lem:moments-t-k}:} First we prove part (i). Recall $\hat t_k$ from \eqref{eqn:def-hat-t-k}. Since $T_k(\gamma_n) \leq |\cE_k| \hat t_k$ and $|\cE_k| \leq 48 \cdot 4^k$ for $k \geq -1$, we have
\begin{equation}
	\label{eqn:1016}
 	\E[T_k^r(\gamma_{n})] \leq \E[ (|\cE_k| \hat t_k)^r] \leq (48 \cdot 4^k)^{r} \E[\hat t_k^r].
\end{equation}
Since for any $r < \eta_0$, there exists $\eta \in (r, \eta_0)$ such that $\E[t_e^{\eta/4}] < \infty$. Recall $c_k$ from Lemma \ref{lem:hat-t-k-moments}. For all $k \geq -1$ we have, $c_k \eta/4 \geq \eta > r$. Thus by Lemma \ref{lem:hat-t-k-moments}, $\E[\hat t_k^r] < \infty$ and (i) is proved.
  
 

Next, we prove (ii). The constants $\epsilon_0, \delta_0$ are from \eqref{eqn:466}. We will perform a decomposition for $\hat p_k$ introduced by J\'arai \cite[p.~319]{jarai2003invasion} using iterated logarithms.  Define $\log^{(0)}k=k$ and $\log^{(j)}k = \log^{(j-1)}k$ for $j \geq 1$ such that it is well-defined. For $k > 10$, let
 \begin{equation*}
 	\log^* k = \min\set{j >0: \log^{(j)}k \mbox{ is well-defined and } \log^{(j)}k \leq 10 }.
 \end{equation*}
 Denote, for $j = 0,1,2,\cdots, \log^* k$,
 \begin{equation*}
 	q_k(j) := p_{\lfloor 2^k/(C_1\log^{(j)}k) \rfloor},
 \end{equation*}
 where $C_1$ is so large that
 \begin{align}
 	&C_1 > 2/\log 10,  \label{eqn:c1-1}\\
  	&2 r \log2  - K_{\ref{eqn:outlet}.2}C_1/2 < -\lambda, \label{eqn:c1-2}\\
 	& \lceil 2r \rceil -K_{\ref{eqn:outlet}.2} C_1 /2 < -1. \label{eqn:c1-3}
 \end{align}
 Given $C_1$, let $k_0 > 10$ be the smallest integer such that for all $k \geq k_0$,
 \begin{align}
 	\label{eqn:def-k0}
 	2^{k/2-1} > C_1 k, \qquad  p_c + 2^{-\epsilon_0  k/2} < 1, \qquad \mbox{ and } k > \frac{\log r}{\log 2} + 3.
 \end{align}
 The reason for the above choices will be clear as the proof proceeds.  We assume $k \geq k_0$ for the rest of the proof. By Lemma \ref{lem:hat-t-k-moments}, the third condition in the above display implies that $\E[\hat t_k^{2r}] \leq K_{\ref{lem:hat-t-k-moments}}(2r,1/4,F)$, and combining this with \eqref{eqn:1016} we have for all $k \geq k_0$,
 \begin{equation}
 	\label{eqn:clever-moment-bound}
 	\E[T_k^{2r}(\gamma_{n})] \leq (48 \cdot 4^k)^{2r} K_{\ref{lem:hat-t-k-moments}}.
 \end{equation}
 Note that $q_k(\log^* k) < \cdots < q_k(1)$ are well-defined as long as $2^k > C_1 k$. We decompose $\E[T_k^r(\gamma_n)]$ as follows:
 \begin{align}
 	\E[T_k^r(\gamma_n)] =& \E\bbrac{T_k^r(\gamma_n) \ind\set{\hat p_k > q_k(0)}} + \sum_{j=0}^{\log^* k -1} \E\bbrac{T_k^r(\gamma_n) \ind\set{q_k(j+1) < \hat p_k \leq q_k(j)}} \nonumber \\
 	&+ \E\bbrac{T_k^r(\gamma_n) \ind\set{\hat p_k \leq q_k(\log^* k)}}. \label{eqn:decomposition}
 \end{align}
 By \eqref{eqn:jarai-correlation-length} and the fact that $C_1 > 2/\log 10$, for $j=0,1,\ldots, \log^* k$ and $k \geq k_0$,
 \begin{equation*}
 	L(q_k(j))  \leq \left\lfloor  \frac{2^k}{C_1\log^{(j)}k} \right\rfloor \leq \frac{2^k}{C_1\log^{(\log^* k)}k} \leq \frac{2^k}{C_1 \log 10} < 2^{k-1}.
 \end{equation*}
 Then applying Lemma \ref{lem:main-lemma} and Lemma \ref{lem:498}, for all $\alpha \geq 1$, $k_0 \leq k \leq n-1$ and $j=0,1,\ldots, \log^* k$,
 \begin{align}
 	\E\bbrac{T_k^{\alpha}(\gamma_n) \ind\set{\hat p_k \leq q_k(j)} } \leq& 
 	 [F^{-1}(q_k(j))]^\alpha \E\bbrac{N^\alpha(2^{k-1}, 2^k, q_k(j))}  \nonumber\\
 	\leq&  [F^{-1}(q_k(j))]^\alpha \cdot \lceil \alpha \rceil! \bpare{\frac{K_{\ref{lem:498}} 2^k}{L(q_k(j))}}^{2\lceil \alpha \rceil} .\label{eqn:560}
 \end{align}
 By \eqref{eqn:def-k0} we have for $k \geq k_0$ and $j=0,\ldots, \log^* k$,
 \begin{equation}
 	\label{eqn:564-1}
  	 \left\lfloor \frac{2^k}{C_1\log^{(j)}k} \right\rfloor \geq \frac{2^{k-1}}{C_1\log^{(j)}k} \geq \frac{2^{k-1}}{C_1 k } > 2^{k/2}.
  \end{equation}
  Then by \eqref{eqn:466} we have
  \begin{equation}
 	 \label{eqn:564-2}
  	q_k(j) \leq p_c + \left\lfloor \frac{2^k}{C_1\log^{(j)}k}  \right\rfloor^{-\epsilon_0 } \leq p_c + 2^{-k\epsilon_0 /2} < 1.
  \end{equation}
  Applying \eqref{eqn:jarai-correlation-length} and \eqref{eqn:564-2} in  \eqref{eqn:560}, we have for $k_0 \leq k \leq n-1$ and $j=0, \ldots, \log^* k$ with $\alpha \geq 1$,
  \begin{equation}
 	 \label{eqn:key-step}
  	\E\bbrac{T_k^{\alpha}(\gamma_n) \ind\set{\hat p_k \leq q_k(j)} } \leq \lceil \alpha \rceil! (C_2 \log^{(j)}k)^{2\lceil \alpha \rceil} a_k^\alpha,
  \end{equation}
  where $C_2 := 2K_{\ref{lem:498}}C_1/K_{\ref{eqn:jarai-correlation-length}} $.
  Now we bound the sum in \eqref{eqn:decomposition}, starting with the last term. Applying \eqref{eqn:key-step} with $\alpha = r$ and $j = \log^* k$, we have for $k_0 \leq k \leq n-1$ and $r \geq 1$,
  \begin{equation}
 	 \label{eqn:581-3}
  	\E\bbrac{T_k^{r}(\gamma_n) \ind\set{\hat p_k \leq q_k(\log^* k)} } \leq \lceil r \rceil! (10C_2)^{2\lceil r\rceil} a_k^r.
  \end{equation}
  For the first term in \eqref{eqn:decomposition}, applying Cauchy-Schwarz inequality, \eqref{eqn:clever-moment-bound} and \eqref{eqn:outlet}, for $k_0 \leq k \leq n-1$,
  \begin{align}
  	\E\bbrac{T_k^{r}(\gamma_n) \ind\set{\hat p_k > q_k(0)} } \leq& 
 	 \E\bbrac{T_k^{2r}(\gamma_n)}^{1/2} \bbrac{\pr(\hat p_k > q_k(0))}^{1/2} \nonumber\\
 	 \leq&  \bpare{(48 \cdot 4^{k})^{2r} K_{\ref{lem:hat-t-k-moments}}}^{1/2} \cdot K_{\ref{eqn:outlet}.1}^{1/2} \exp(-K_{\ref{eqn:outlet}.2}C_1 k/2) \nonumber\\
 	 =&  48^r \bpare{K_{\ref{lem:hat-t-k-moments}} K_{\ref{eqn:outlet}.1}}^{1/2}  \exp\bpare{2 r k\log2  - K_{\ref{eqn:outlet}.2}C_1 k/2}. \label{eqn:581-1}
  \end{align}
  For the second term in \eqref{eqn:decomposition}, applying Cauchy-Schwarz inequality, \eqref{eqn:key-step} with $\alpha =2r$, and \eqref{eqn:outlet}, we have for $j=0,1,\cdots,\log^* k -1$ and $k_0 \leq k \leq n-1$,
  \begin{align}
  	&\E\bbrac{T_k^r(\gamma_n) \ind\set{q_k(j+1) < \hat p_k \leq q_k(j)}} \nonumber\\
 	\leq&  \E\bbrac{T_k^{2r}(\gamma_n) \ind\set{\hat p_k \leq q_k(j)}}^{1/2} \bbrac{\pr(\hat p_k > q_k(j+1))}^{1/2} \nonumber\\
 	\leq&  \bbrac{\lceil 2r \rceil! (C_2 \log^{(j)}k)^{2\lceil 2r \rceil} a_k^{2r}}^{1/2}\cdot K_{\ref{eqn:outlet}.1}^{1/2} \exp(-K_{\ref{eqn:outlet}.2} C_1  \log^{(j+1)} k /2) \nonumber\\
 	=& (\lceil 2r \rceil !)^{1/2} C_2^{ \lceil 2r \rceil} K_{\ref{eqn:outlet}.1}^{1/2} a_k^r (\log^{(j)}k)^{\lceil 2r \rceil -K_{\ref{eqn:outlet}.2} C_1 /2}. \label{eqn:581-2}
  \end{align}
 Then combining \eqref{eqn:581-1}, \eqref{eqn:581-2}, \eqref{eqn:581-3} and using the definition of $C_1$ in \eqref{eqn:c1-2} and \eqref{eqn:c1-3}, there are $C_3, C_4, C_5 >0$ such that for $ k_0 \leq k \leq n-1$,
 \begin{equation*}
 	\E[T_k^r(\gamma_n)] \leq C_3 e^{- \lambda k} + C_4 a_k^r  \sum_{j=0}^{\log^* k -1} (\log^{(j)} k)^{-1} + C_5 a_k^r.
 \end{equation*}
 Using \cite[Eq. 2.16]{van2007size} which says $\sum_{j=0}^{\log^* k} (\log^{(j)}k)^{-1}$ is uniformly bounded in $k$, we complete the proof of Lemma \ref{lem:moments-t-k}. \qed

\section{Study of the mean}
\label{sec:proof-mean}

In this section, we give the proof of Theorem \ref{thm:mean}. We prove Corollary \ref{cor:rho-f} in Section~\ref{sec:proof-rho-f}.

\subsection{Proof of Theorem \ref{thm:mean}}
\label{sec:proof-thm-mean}

First we prove an elementary lemma.

\begin{lem}\label{lem:sandwich-series}
	Let $f(t)$, $t \in [0,\infty)$, be a positive non-increasing function. Fix $\delta > \epsilon >0$ and integers $k_1, k_2 \geq 1$. Then there exist constants $C_1, C_2 \in (0,\infty)$ such that for all $n \geq \max\set{k_1,k_2}$,
	\begin{equation*}
		 C_1 \sum_{k=k_2}^n f(\epsilon k) \leq \sum_{k=k_1}^n f(\delta k) \leq C_2 \sum_{k=k_2}^n f(\epsilon k).
	\end{equation*}
\end{lem}
\textbf{Proof: } It suffices to prove the lemma for $k_1 = k_2 = 1$. If $ \delta k \leq \epsilon  k' < \delta (k + 1)$, then $f(\epsilon  k' ) \leq f(\delta k)$, and for each $k$ there are at most $\lceil  \delta / \epsilon  \rceil$ such integer $k'$. Therefore
\[
\sum_{k=1}^n f(\epsilon k) \leq \lceil \delta / \epsilon \rceil (f(\epsilon) + \sum_{k=1}^n f(\delta k)).
\]
As $\epsilon < \delta$, we obtain

	\begin{equation*}
		1 \leq \frac{\sum_{k=1}^{n} f(\epsilon  k)}{ \sum_{ k = 1}^n f(\delta k)} \leq \frac{\lceil \delta / \epsilon  \rceil \left(f(\epsilon) +  \sum_{ k = 1}^n f(\delta k)\right)}{ \sum_{ k = 1}^n f(\delta k)},
	\end{equation*}
	and this completes the proof.
\qed\\

 Now we prove Theorem \ref{thm:mean}.

\textbf{Proof of Theorem \ref{thm:mean}:} For the upper bound, note that for $n \geq -1$,
\begin{equation*}
	\E T(\vzero,\partial B(2^n)) \leq  \sum_{k=-1}^{n-1} \E T_{k}(\gamma_n) .
\end{equation*}
Take $k_0$ as in Lemma \ref{lem:moments-t-k} and apply this lemma with $r = 1$. We obtain $\E T(\vzero,\partial B(2^n)) < \infty$ for all $n \geq -1$ and in particular for $n \geq k_0 + 1$,
\begin{align*}
\E T(\vzero,\partial B(2^n)) \leq \sum_{k=-1}^{k_0-1}  \E T_{k}(\gamma_n) + K_{\ref{lem:moments-t-k}} \sum_{k=k_0}^{n-1} [a_k + e^{-k}] \leq C_1 \sum_{k=2}^n F^{-1}(p_c + 2^{-k}),
\end{align*}
for some constant $C_1 >0$. The last inequality uses Lemma \ref{lem:sandwich-series} and the fact $F^{-1}(p_c+1/4) > 0$. This proves (i). 

For the lower bound, the proof is similar to that of \cite[Theorem~8.1.2]{zhang1999double}. 
By \eqref{eqn:466}, for $k\geq 1$, crossing a $p_{2^k}$-closed dual circuit costs passage time at least $F^{-1}(p_c+2^{-k\delta_0})$. Therefore, by \eqref{eqn:5-26}, $\E T(\vzero,\partial B(2^n))$ is bounded below by
\begin{align*}
& \sum_{k=1}^n \pr(\text{there is a }p_{2^k}\text{-closed dual circuit around }\vzero \text{ in }B(2^k)^*\setminus B(2^{k-1})^*) \cdot F^{-1}(p_c+2^{-k\delta_0}) \\
\geq ~&\sum_{k=1}^n K_{\ref{eqn:5-26}}F^{-1}(p_c+2^{-k\delta_0}) .
\end{align*}
Applying Lemma \ref{lem:sandwich-series} completes the proof of (ii). \qed

\subsection{Proof of Corollary~\ref{cor:rho-f}}
\label{sec:proof-rho-f}

To prove Corollary~\ref{cor:rho-f}, we need the following definition from \cite[p.~146]{zhang1999double}. Given two distribution functions $G$ and $H$, we say that $G\preceq H$ if there exists $\xi>0$ such that $G(x)\leq H(x)$ for all $0\leq x\leq \xi$. By \cite[Theorem~8.1.4]{zhang1999double}, if $\rho(G)<\infty$ almost surely and if $G\preceq H$, then $\rho(H)<\infty$ almost surely. (This is in fact provable in general dimensions, though Zhang only gave a proof for $d=2$.)

\textbf{Proof of Corollary~\ref{cor:rho-f}: }
Suppose that $\sum_{n=2}^\infty F^{-1}(p_c+2^{-n})<\infty$. Let $\xi>0$ be arbitrary, and let $\tilde{F}$ be a distribution function such that $\tilde{F}=F$ on $[0,\xi]$ and there exists $x_0>0$ such that $\tilde{F}(x_0)=1$. Note that we still have $\sum_{n=2}^\infty \tilde{F}^{-1}(p_c+2^{-n})<\infty$ and $\tilde{F}$ has all moments. By Theorem~\ref{thm:mean}(i), $\rho(\tilde{F})<\infty$ almost surely. Since $\rho(F)\leq \rho(\tilde{F})$, we have $\rho(F)<\infty$ almost surely.

Now suppose that $\sum_{n=2}^\infty F^{-1}(p_c+2^{-n})=\infty$. This implies that $\sum_k F^{-1}(p_c + 2^{-\delta_0 k}) = \infty$ for $\delta_0$ from \eqref{eqn:466}. For $k\geq 1$, write
\[
A_k = \{\exists\;p_{2^k}\text{-closed dual circuit around $\vzero$ in }B(2^k)^*\setminus B(2^{k-1})^*\}\ 
\]
and $b_k := F^{-1}(p_c + 2^{-k\delta_0})$. For $n \geq 1$, define $S_n = \sum_{k=1}^n b_k \ind_{A_k}$ and compute by Cauchy-Schwarz, \eqref{eqn:5-26} and independence of the $A_k$'s:
\[
\E S_n^2 \leq \sum_{j,k=1}^n b_jb_k = \left( \sum_{k=1}^n b_k \right)^2 \leq \frac{1}{K_{\ref{eqn:5-26}}^2} \left(\E S_n \right)^2\ .
\]
By the Paley-Zygmund inequality (second moment method), we can find $D>0$ such that for all $n \geq 1$,
\[
\pr\left( S_n \geq D \E S_n \right) > D\ .
\]
Therefore
\[
\pr\left( S_n \geq D \E S_n \text{ for infinitely many }n \right) \geq D\ .
\]
Because $\rho(F) \geq S_n$ for all $n \geq 1$, this implies that $\rho(F) = \infty$ with probability at least $D$. Since, $\rho(F) < \infty$ almost surely or $\rho(F) = \infty$ almost surely, this completes the proof. \qed\\

\section{Study of the variance}
\label{sec:proof-variance}

In this section, we prove Theorems \ref{thm:variance} and \ref{thm:limit-theorem}. The main tool is the construction of a martingale introduced in \cite{kesten1997central}. We start with some definitions.

 Define the annuli
\begin{equation*}
	\Ann(n)  = B(2^{n+1}) \setminus B(2^n), \mbox{ for } n \geq 0 \text{ and } \Ann(-1) = B(1).
\end{equation*}
For a vertex self-avoiding circuit $\cC$ in $\bZ^2$, write $\bar \cC$ for the graph induced by all the vertices in $\bZ^2$ that are either on or in the interior of $\cC$. Define
\begin{equation}
	\label{eqn:def-m}
	m(n) := \inf\set{ k \geq n: \text{There is a $p_c$-open circuit in $\Ann(k)$ surrounding $\vzero$} }, \mbox{ for } n \geq -1.
\end{equation}
Note that $m(n) \geq n$. Write $m(n)=m(n,\omega)$ to emphasize the underlying weights $\omega \in \Omega$.
Denote $\cC_n$ as the smallest $p_c$-open circuit in $\Ann(m(n))$: precisely,
\begin{equation}
	\label{eqn:def-c-n}
	\cC_n := \mbox{ the innermost $p_c$-open circuit } \cC \subset \Ann(m(n)) \mbox{ surrounding $\vzero$ }.
\end{equation}
Define
\begin{equation}
	\label{eqn:fn-filtration}
	\cF_n := \text{sigma-field generated by $\cC_n$ and $\set{\omega_e: e \in \bar \cC_n}$}.
\end{equation}
By definition we have $\cC_{n}(\omega) = \cC_{m(n,\omega)}(\omega)$. For $n < n'$, we have $m(n) \leq m(n')$, thus $\set{\cF_n}_{n \in \bN}$ forms a filtration. Denote $\cF_{-1} = \set{\emptyset, \Omega}$ and $\cC_{-1} = \{\vzero\}$. Instead of $T(\vzero, \partial B(2^n))$, we first try to study $T(\vzero, \cC_n)$. Write
\begin{align*}
	T(\vzero,\cC_n) - \E T(\vzero, \cC_n) 
	= \sum_{ k = 0}^n\left( \E[ T(\vzero, \cC_n) \mid \cF_k] - \E[ T(\vzero, \cC_n) \mid \cF_{k-1}]\right) =: \sum_{k=0}^n \Delta_k.
\end{align*}
Then $\set{\Delta_k}_{0 \leq k \leq n}$ is an $\cF_k$-martingale increment sequence. Thus
\begin{equation}
	\label{eqn:1061}
	\var(T(\vzero,\cC_n)) = \sum_{k=0}^n \E[\Delta_k^2].
\end{equation}

The following two theorems are the results for $T(\vzero, \cC_n)$ corresponding to those in Theorems~\ref{thm:variance} and \ref{thm:limit-theorem}.
\begin{thm} 
	\label{thm:variance-zero-circuit}
	Let $\eta_0$ be as defined in \eqref{eqn:def-eta0}. \\
	(i) If $\eta_0 >2$, then there exists $C_1 > 0$ such that for $n \geq 2$,
	$$ \var(T(\vzero,\cC_n)) \leq C_1 \sum_{k=2}^n [F^{-1}(p_c + 2^{-k})]^2.$$
	(ii) There exists $C_2 > 0$ such that for $n \geq 2$,
	$$\var(T(\vzero,\cC_n)) \geq C_2 \sum_{k=2}^n [F^{-1}(p_c + 2^{-k})]^2.$$
	
\end{thm}

\begin{thm}
	\label{thm:clt-zero-circuit}
	Assume that $\eta_0 > 2$. Further assume $\sum_{k = 1 }^\infty [F^{-1}(p_c + 2^{-k})]^2 = \infty$. Then as $n \to \infty$,
		\begin{equation*}
		\frac{T(\vzero, \cC_n) - \E T(\vzero, \cC_n)}{(\var(T(\vzero, \cC_n)))^{1/2}}
		= \frac{\sum_{k=0}^n \Delta_{k}}{(\sum_{k=0}^n \E[\Delta_{k}^2])^{1/2}} \convd N(0,1).
	\end{equation*}
\end{thm}

We will first prove Theorem \ref{thm:variance-zero-circuit} in Section \ref{sec:proof-variance-zero-circuit}. Next in Section \ref{sec:proof-clt-zero-circuit} we prove the CLT stated in Theorem \ref{thm:clt-zero-circuit}. Finally, in Section \ref{sec:proof-variance-clt-box}, we control the difference between $T(\vzero, \cC_q)$ and $T(\vzero, \partial B(n))$ for $2^{q-1} \leq n \leq 2^q-1$ and prove Theorems \ref{thm:variance} and \ref{thm:limit-theorem}.

\subsection{Proof of Theorem \ref{thm:variance-zero-circuit}}
\label{sec:proof-variance-zero-circuit}

Because of \eqref{eqn:1061}, we first study bounds on the moments of $\Delta_k$. An important ingredient is an alternative formula for $\Delta_k$ which was proved in \cite[Lemma~2]{kesten1997central}, and we state it as part (ii) in the following lemma. Denote $(\Omega', \cF', \pr')$ as another copy of the probability space $(\Omega, \cF, \pr)$. Let $\E'$ denote the expectation with respect to $\pr'$, and $\omega'$ denote a sample point in $\Omega'$. Denote $m(n,\omega)$, $\cC_k(\omega)$ and $T(\cdot,\cdot)(\omega)$ for the quantities defined as in the previous sections, but with explicit dependence on $\omega$. Define $\ell(n,\omega, \omega') := m(m(n,\omega)+1, \omega')$, or equivalently,
\begin{equation*}
	\ell(n,\omega,\omega') := \inf \set{ \ell \geq m(n,\omega) + 1: \exists \text{ $p_c$-open circuit around } \vzero \text{ in $\Ann(\ell)$ in the weights $\omega'$ }}.
\end{equation*}
We need the following results, which are \cite[Lemma~3]{kesten1997central} and \cite[Lemma~2]{kesten1997central}.
\begin{lem}[Kesten-Zhang]
	\label{lem:delta-p-decomposition}
	(i)  There exists $K_{\ref{lem:delta-p-decomposition}} > 0$ such that for all integers $k,t \geq 1$,
	 $$ \pr(m(k) \geq k + t) \leq \exp(-K_{\ref{lem:delta-p-decomposition}} t).  $$
	(ii) For all $k \geq 0$, $\Delta_k$ does not depend on $n$. Precisely,
	\begin{align*}
		\Delta_k(\omega) = T(\cC_{k-1}(\omega), \cC_{k}(\omega)) + \E'[T(\cC_k(\omega), \cC_{\ell(k,\omega,\omega')}(\omega'))(\omega') ] - \E'[T(\cC_{k-1}(\omega), \cC_{\ell(k,\omega,\omega')}(\omega')) (\omega')].
	\end{align*}
	
\end{lem}
Sometimes we will write $T(\cdot,\cdot)$ instead of $T(\cdot,\cdot)(\omega)$ or $T(\cdot,\cdot)(\omega')$ when the meaning is clear from the context. The following lemma is a consequence of Lemma \ref{lem:moments-t-k}.
\begin{lem} 
	\label{lem:483}
	Assume $\eta_0 > 1$. 
	
	(i) For any $r \in [1,\infty)$ and $\lambda \in (0,\infty)$, there exist $k_0 = k_0 (r,\lambda, F)>0$ and $K_{\ref{lem:483}} = K_{\ref{lem:483}}(r) =   K_{\ref{lem:483}}(r,\lambda, F)> 0$ such that for all $k \geq k_0$ and $\ell \geq 1$,
	\begin{equation*}
		\E \left[ T^r( \partial B(2^k),  \partial B(2^{k+\ell}) )  \right] \leq K_{\ref{lem:483}} \ell^{r} (a_k^r + e^{-\lambda k}).
	\end{equation*}
	
	(ii) For any $r \in [1,\eta_0)$, there exists a constant $K_{\ref{lem:483}} = K_{\ref{lem:483}}(r) = K_{\ref{lem:483}}(r, F) > 0$ such that for all $k \geq -1$ and $\ell \geq 1$, 
	\begin{equation*}
		\E \left[ T^r( \partial B(2^k),  \partial B(2^{k+\ell}) )  \right] \leq K_{\ref{lem:483}} \ell^{r}.
	\end{equation*}
\end{lem}
\textbf{Proof: } Take $n := k+\ell + 1$. Since $\gamma_n \cap (B(2^{k + \ell}) \setminus B(2^k))$ provides a specific path connecting the inner and outer boundaries of the annulus, we have $T( \partial B(2^k),  \partial B(2^{k+\ell}) ) \leq \sum_{i = k}^{k + \ell -1} T_i(\gamma_n)$. Therefore, by applying Jensen's inequality and Lemma \ref{lem:moments-t-k} for $k \geq k_0$, we have
\begin{equation*}
	\E \left[ T^r( \partial B(2^k),  \partial B(2^{k+\ell}) )  \right] \leq \ell^{r-1} \left(\sum_{i = k}^{k+\ell - 1} \E[T_i^r(\gamma_n)] \right) \leq \ell^{r-1} \cdot K_{\ref{lem:moments-t-k}} \sum_{i=k}^{k+\ell -1} (a_i^r + e^{-\lambda i} )\leq K_{\ref{lem:moments-t-k}} \ell^{r} (a_k^r + e^{-\lambda k}).
\end{equation*}
This proves (i). To prove (ii), if $r \in [0,\eta_0)$, by Lemma \ref{lem:moments-t-k} (parts (i) and (ii) combined), for all $n \geq k \geq -1$ we have $\E[T_k^r(\gamma_n)] < C_1$ for some constant $C_1 = C_1(r,F)>0$. Using this fact in the above bound proves (ii).  \qed\\

We now show how the above lemma implies bounds on moments of the $\Delta_k$'s.
\begin{lem}
	\label{lem:delta-p-upper-bound}
	Assume $\eta_0 > 1$. 
	
	(i) For any $r \in [1, \infty)$ and $\lambda \in (0,\infty)$, there 
	exist $K_{\ref{lem:delta-p-upper-bound}} = K_{\ref{lem:delta-p-upper-bound}}(r) = K_{\ref{lem:delta-p-upper-bound}}(r, \lambda , F) > 0$ and $k_0 = k_0(r,\lambda, F) > 0$ such that for all  $k \geq k_0 + 1$,
	\begin{equation*}
		\E[|\Delta_k|^r] \leq K_{\ref{lem:delta-p-upper-bound}} \cdot (a_{k-1}^r + e^{-\lambda k}).
	\end{equation*}
	
	(ii) For any $r \in [1,\eta_0)$ and $k \geq 0$, we have $\E[|\Delta_k|^r] < \infty$.
\end{lem}

\textbf{Proof: } First we prove (i).
Using  the fact that $T(\cC_{k}(\omega), \cC_{\ell(k,\omega,\omega')}(\omega'))  \leq T(\cC_{k-1}(\omega), \cC_{\ell(k,\omega,\omega')}(\omega')) $ and Lemma \ref{lem:delta-p-decomposition} (ii), we have
\begin{equation*}
	|\Delta_k(\omega)| \leq T(\cC_{k-1}(\omega), \cC_{k}(\omega)) + \E'[T(\cC_{k-1}(\omega), \cC_{\ell(k,\omega,\omega')}(\omega')) ].
\end{equation*}
By Jensen's inequality, we have

\begin{equation}
	\label{eqn:503}
	\frac{1}{2^{r-1}}\E|\Delta_k(\omega)|^r \leq  \E \bbrac{T^r(\cC_{k-1}(\omega), \cC_{k}(\omega))} + \E \bbrac{\left(\E'[T(\cC_{k-1}(\omega), \cC_{\ell(k,\omega,\omega')}(\omega')) ]\right)^r}.
\end{equation}
First we give an upper bound for the second term.  Recall $k_0(r,\lambda, F)$ from Lemma \ref{lem:delta-p-upper-bound}. Fix $\omega \in \Omega$, and estimate for $k \geq k_0(2,\lambda, F)+1$,
\begin{align}
	&~\E'[T(\cC_{k-1}(\omega), \cC_{\ell(k,\omega,\omega')}(\omega')) ] \nonumber\\
	=&~ \sum_{t = 0}^\infty \E'\bbrac{ T(\cC_{k-1}(\omega), \cC_{\ell(k,\omega,\omega')}(\omega'))  \ind_{\set{ \ell(k,\omega, \omega') - m(k,\omega) -1 = t}} }  \nonumber\\
	\leq&~ \sum_{t = 0}^\infty \E'\bbrac{ T(\partial B(2^{k-1}), \partial B(2^{m(k,\omega) +2 + t}))\ind_{\set{ \ell(k,\omega, \omega') - m(k,\omega) -1 = t}}}  \nonumber\\
	\leq&~ \sum_{t = 0}^\infty \E'\bbrac{ T^2(\partial B(2^{k-1}), \partial B(2^{m(k,\omega) +2 + t}))}^{1/2} \pr' \bpare{ \ell(k,\omega, \omega') - m(k,\omega) -1 = t}^{1/2}  \nonumber\\
	\leq&~ \sum_{t = 0}^\infty (K_{\ref{lem:483}}(2))^{1/2} (a_{k-1}^2 + e^{-\lambda (k-1)})^{1/2} (m(k,\omega)-k + t +3) e^{-{K_{\ref{lem:delta-p-decomposition}}t}/2} \nonumber\\
	\leq&~ C_1 (m(k,\omega)-k + 1 )(a_{k-1} + e^{-\lambda (k-1)/2}), \label{eqn:691}
\end{align}
where the fourth line uses the Cauchy-Schwarz inequality, the fifth line uses Lemma \ref{lem:483} with $r = 2$ and Lemma \ref{lem:delta-p-decomposition}(i), and in the fifth line $C_1 := {K_{\ref{lem:483}}(2)}^{1/2} \sum_{t=0}^\infty (t+2)e^{-K_{\ref{lem:delta-p-decomposition}}t/2}$. Therefore
\begin{align}
	\E \bbrac{\left(\E'[T(\cC_{k-1}(\omega), \cC_{\ell(k,\omega,\omega')}(\omega')) ]\right)^r}  
	\leq& C_1^r (a_{k-1} + e^{-\lambda (k-1)/2})^r \E[(m(k,\omega)-k + 1 )^r]  \nonumber\\
	\leq& C_1^r \E[(m(k,\omega)-k + 1 )^r]\cdot 2^{r-1}(a_{k-1}^r + e^{-\lambda r (k-1)/2} ). \label{eqn:698}
\end{align}
By Lemma \ref{lem:delta-p-decomposition}(i) $\E[(m(k,\omega)-k + 1 )^r] < \infty$, 
so we complete the bound on the second term in \eqref{eqn:503}. To bound the 
first term in \eqref{eqn:503}, similar to \eqref{eqn:691}, applying 
Cauchy-Schwarz inequality, we have for $k \geq k_0(2r,\lambda, F)+1$,
\begin{align}
	 \E \bbrac{T^r(\cC_{k-1}(\omega), \cC_{k}(\omega))} 
	 \leq& \sum_{t=0}^\infty \E[T^{2r}(\partial B(2^{k-1}), \partial B(2^{k + t + 1}))]^{1/2} \pr \bpare{m(k)-k = t}^{1/2} \nonumber\\
	 \leq& \sum_{t=0}^\infty [K_{\ref{lem:483}}(2r)]^{1/2} (a_{k-1}^{2r} + e^{-\lambda(k-1)})^{1/2} (t+2)^r \cdot e^{-{K_{\ref{lem:delta-p-decomposition}}t}/{2}} \nonumber\\
	 \leq& (a_{k-1}^r + e^{-\lambda (k-1)/2}) \bpare{[K_{\ref{lem:483}}(2r)]^{1/2}\sum_{t=0}^\infty (t+2)^r e^{-{K_{\ref{lem:delta-p-decomposition}}t}/2}}. \label{eqn:707}
\end{align}
Combining \eqref{eqn:503}, \eqref{eqn:698} and \eqref{eqn:707}, we complete the proof of Lemma \ref{lem:delta-p-upper-bound} (i). 
The proof of part (ii) can be done in exactly the same way, using Lemma \ref{lem:483} (ii).  \qed\\

\begin{figure}[t]
\begin{center}

\includegraphics[trim=1.5cm 7cm 1.5cm 3cm, clip=true, width=10cm]{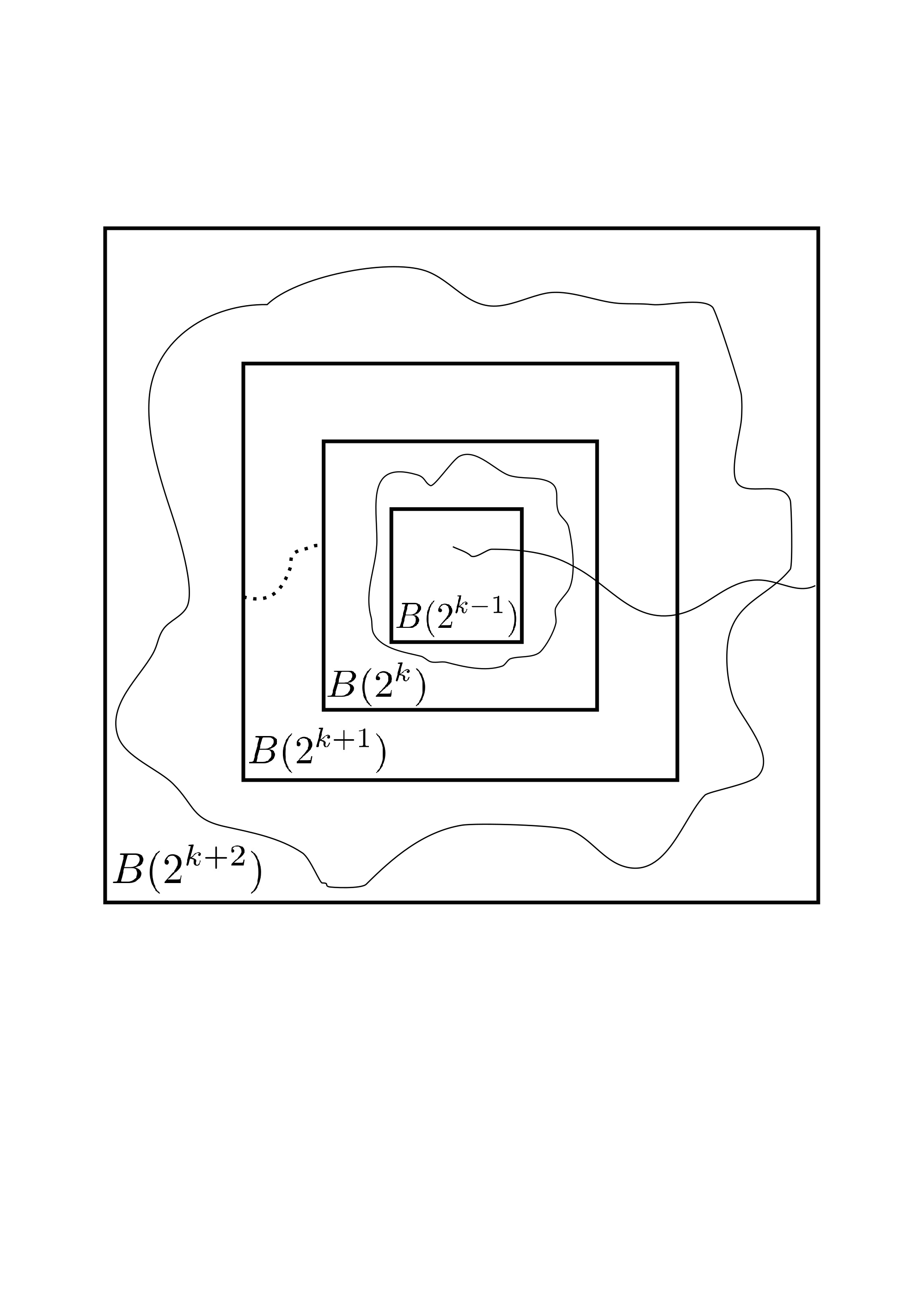}
\caption{The events (1)-(4) in the proof of Lemma~\ref{lem:delta-p-lower-bound}. The $p_c$-open crossing on the right connects the two $p_c$-open circuits around the origin, but the $p_c$-closed path on the left (shown as a dotted curve) blocks the existence of a $p_c$-open circuit around $\vzero$ in $B(2^{k+1}) \setminus B(2^k)$.}
\label{fig: lower_bound}
\end{center}
\end{figure}

The next lemma gives a lower bound for $\E[\Delta_k^2]$.
\begin{lem}
	\label{lem:delta-p-lower-bound}
	There exists a constant $K_{\ref{lem:delta-p-lower-bound}} > 0$ such that for all integers $k \geq 2 $, we have
	\begin{equation*}
		\E[\Delta_k^2] \geq K_{\ref{lem:delta-p-lower-bound}} [F^{-1}(p_c+ 2^{-\delta_0 k})]^2,
	\end{equation*}
	where $\delta_0$ is from \eqref{eqn:466}.
\end{lem}
\textbf{Proof: } Recall the expression of $\Delta_k$ in Lemma \ref{lem:delta-p-decomposition}(ii) and the filtration $\cF_k$ in \eqref{eqn:fn-filtration}. The goal of the proof is to construct an event $E \in \cF_k$ with $\pr(E) >0$ uniformly in $n,k$ such that for $\omega \in E$,
\begin{align}
	&T(\cC_{k-1}(\omega), \cC_k(\omega))(\omega) = 0, \\
	&\E'[T(\cC_{k-1}(\omega), \cC_{\ell(k,\omega,\omega')}(\omega')) ] - \E'[T(\cC_k(\omega), \cC_{\ell(k,\omega,\omega')}(\omega'))] > C_2 F^{-1}(p_c + 2^{-\delta_0 k}), \label{eqn:578}
\end{align}
where $C_2 > 0$ is a constant. We start by defining the event $\tilde E$ to be the intersection of the following events (see Figure~\ref{fig: lower_bound}):
\begin{enumerate}
	\item There exists a $p_c$-open circuit around $\vzero$ in $B(2^{k})\setminus B(2^{k-1})$.
	\item There exists a $p_c$-open circuit around $\vzero$ in $B(2^{k+2})\setminus B(2^{k+1})$.
	\item There exists a $p_c$-open left-right crossing of $[0, 2^{k+2}] \times [-2^{k-1}, 2^{k-1}]$.
	\item There exists a \emph{dual} $p_c$-closed left-right crossing of $[-2^{k+1}, -2^{k}]^* \times [-2^{k}, 2^{k}]^*$.
\end{enumerate}

By the RSW Theorem (\cite[Sect. 11.7]{grimmett1999percolation}), each of the above events has probability bounded from below for all $k \geq 1$. The events (1), (2) and (3) are all non-increasing, and they are jointly independent from (4). Therefore applying independence and the FKG inequality, there exists a constant $C_3 > 0$ such that $\pr(\tilde E) \geq C_3$ for all $k\geq 1$ . Now consider a new event $(3')$:
\begin{enumerate}
	\item[$(3')$] There exists a $p_c$-open left-right crossing of $\bar \cC_k \cap ([0, 2^{k+2}] \times [-2^{k-1}, 2^{k-1}])$.
\end{enumerate}
Define the event $E$ to be the intersection of the events (1),(2), (4) and $(3')$. Then $E \in \cF_k$, $\tilde E \subset E$ and therefore $\pr(E) \geq \pr(\tilde E) \geq C_3 > 0$. By definition, we have $T(\cC_{k-1}(\omega), \cC_k(\omega))(\omega) = 0 $ for $ \omega \in E$. To see \eqref{eqn:578}, recall the definition of $p_n$ in \eqref{eqn:def-pn}. Let $E' \subset \Omega'$ be the event
\begin{equation*}
	E' = \set{ \text{ There exists a $p_{2^k}$-closed dual circuit around $\vzero$ in $B(2^{k+1})^* \setminus B(2^{k})^*$ }}.
\end{equation*}
From \eqref{eqn:5-26}, let $C_2 > 0$ be such that $\pr'(E') > C_2$ for all $k$. When $\omega \in E$ and $\omega' \in E'$, since every path between $\cC_{k-1}(\omega)$ and $\cC_{\ell(k,\omega,\omega')}(\omega')$ must cross the $p_{2^k}$-closed dual circuit defined in $E'$ and then cross $\cC_k(\omega)$, we have for $k \geq 2$
\begin{equation*}
	T(\cC_{k-1}(\omega), \cC_{\ell(k,\omega,\omega')}(\omega'))(\omega') - T(\cC_{k}(\omega), \cC_{\ell(k,\omega,\omega')}(\omega'))(\omega') \geq F^{-1}(p_{2^k}) \geq F^{-1}(p_c + 2^{-\delta_0 k}),
\end{equation*}
where the last inequality follows from \eqref{eqn:466}. This proves \eqref{eqn:578} and therefore we have $\pr(\Delta_k < - C_2 F^{-1}(p_c + 2^{-\delta_0 k})) \geq \pr(E) \geq C_3$, completing the proof of Lemma \ref{lem:delta-p-lower-bound}. \qed\\

Now we complete the proof of Theorem \ref{thm:variance-zero-circuit}.

\textbf{Proof of Theorem \ref{thm:variance-zero-circuit}: } First we prove (i). Lemma \ref{lem:delta-p-upper-bound} (i) with $r=2$ and $\lambda = 1$ implies that there exists $k_0 \geq 1$ such that for all $k \geq k_0+1$, we have $\E[ \Delta_k^2] \leq K_{\ref{lem:delta-p-upper-bound}} (a_{k-1}^2 + e^{-k})$. For $k \leq k_0$, we will use the general fact that $\E[\Delta_k^2] < C_1$ for some $C_1 > 0$. Therefore for $k \geq k_0 + 1$,
\begin{align*}
	\var(T(\vzero, \cC_k)) 
	= \sum_{k=0}^{k_0} \E[\Delta_k^2] + \sum_{k = k_0+1}^n \E[\Delta_k^2]
	\leq& (k_0+1) C_1 +  K_{\ref{lem:delta-p-upper-bound}} \sum_{k=k_0}^{n-1} (a_k^2 + e^{-k}) \\
	=& C_2 +  K_{\ref{lem:delta-p-upper-bound}} \sum_{k=k_0}^{n-1} a_k^2,
\end{align*}
where $C_2 > 0$ is some constant. Applying Lemma \ref{lem:sandwich-series} with $f(t) := (F^{-1}(p_c + 2^{-t}) \wedge a_{k_0})^2$, $t \geq 0$, completes the proof of (i).

Next we prove (ii). By Lemma \ref{lem:delta-p-lower-bound}, since $\delta_0>1$, we have for $n \geq 2$,
\begin{equation*}
	\var(T(\vzero, \cC_n)) =  \sum_{k=0}^n \E[\Delta_k^2] \geq K_{\ref{lem:delta-p-lower-bound}} \sum_{k=2}^n
	 [F^{-1}(p_c + 2^{-\delta_0 k})]^2 .
\end{equation*}
Applying Lemma \ref{lem:sandwich-series} again completes the proof of (ii).  \qed\\

\subsection{Proof of Theorem \ref{thm:clt-zero-circuit}} 
\label{sec:proof-clt-zero-circuit}

In this section, we prove the CLT for $T(\vzero,\cC_n)$ as stated in Theorem \ref{thm:clt-zero-circuit}. The version of the martingale CLT we use here is from McLeish \cite[Theorem 2.3]{mcleish1974dependent}, which is also used in \cite{kesten1997central}. We state the theorem here.
\begin{thm}[McLeish]
	\label{thm:mart-clt-mcleish}
	Let $\set{X_{k,n}}_{ k \leq n}$ be a martingale difference array satisfying\\
	(i) $\sup_{n} \E[(\max_{k \leq n} X_{k,n})^2] < \infty$,\\
	(ii) $ \max_{k \leq n} |X_{k,n}| \convp 0 $ as $n \to \infty$,\\
	(iii) $\sum_{k \leq n} X_{k,n}^2 \convp 1$ as $n \to \infty$.\\
	Then we have $\sum_{k \leq n} X_{k,n} \convd N(0,1)$.
\end{thm}

We start by introducing some definitions. Using $\set{\Delta_k: k \geq 0}$, define the strong mixing coefficient
\begin{equation*}
	\alpha(\ell) := \sup_{k \geq 0} \sup_{A,B} |\pr(A \cap B) - \pr(A) \pr(B)| \quad \mbox{ for } \ell \geq 0,
\end{equation*}
 where $\sup_{A,B}$ is the supremum over all events $A \in \sigma(\set{\Delta_i: 0 \leq i \leq k})$ and $B \in \sigma(\set{\Delta_i: i \geq k + \ell})$. 

%
%
In order to verify the conditions in Theorem \ref{thm:mart-clt-mcleish}, we first study the strong mixing coefficient $\alpha(\ell)$.
\begin{lem}\label{lem:strong-mixing-coefficient}
	For the process $\set{\Delta_k: k \geq 0}$, there exist $C_1, C_2 > 0$ such that
	\begin{equation*}
		\alpha(\ell) \leq C_1e^{- C_2 \ell} \qquad \mbox{ for } \ell \geq 0.
	\end{equation*}
\end{lem}
\textbf{Proof:} Recall the filtration $\set{\cF_k: k \geq 0}$ associated with $\set{\Delta_k: k \geq 0}$. Denote $\tilde \cF_k := \sigma(\set{\Delta_\ell: \ell \geq k})$ for $k \geq 0$. Now fix $k \geq 0$ and $\ell \geq 2$ with $A \in \cF_{k}$ and $B \in \tilde \cF_{k + \ell}$. By the expression for $\Delta_k$ in Lemma \ref{lem:delta-p-decomposition}(ii), we have
\begin{equation*}
	B \in \sigma\bpare{\set{\omega_e: e \notin E(B(2^{k + \ell - 1}))}}.
\end{equation*}
The next fact was observed in \cite[Lemma 1]{kesten1997central}: 
\begin{equation*}
	A \cap \set{m(k) \leq k + \ell-2} \in \sigma\bpare{\set{\omega_e: e \in E(B(2^{k+\ell-1}))}}.
\end{equation*}
Therefore we have
\begin{equation*}
	\pr(A \cap \set{m(k) \leq k + \ell-2} \cap B) = \pr(A \cap \set{m(k) \leq k + \ell-2}) \pr( B).
\end{equation*}
Thus
\begin{align*}
	|\pr(A \cap B) - \pr(A) \pr(B)| \leq& \pr(A \cap \set{m(k) \geq k + \ell-1}\cap B) + \pr(A \cap \set{m(k) \geq k + \ell-1}) \pr( B) \\
	\leq& 2\pr(m(k) \geq k + \ell-1) \leq 2 \exp( - K_{\ref{lem:delta-p-decomposition}}(\ell - 1)),
\end{align*}
where the last line used Lemma \ref{lem:delta-p-decomposition}(i).
The above bound is also true for $\ell=0,1$ since $|\pr(A \cap B) - \pr(A) \pr(B)| \leq 2$. Also it is independent of $k, A$ and $B$, so this completes the proof for Lemma \ref{lem:strong-mixing-coefficient}. \qed\\

Now we are ready to prove Theorem \ref{thm:clt-zero-circuit}.

\textbf{Proof of Theorem \ref{thm:clt-zero-circuit}: } Since $\eta_0 > 2$, by Lemma \ref{lem:delta-p-upper-bound}, there exist constants $k_1, C_3, C_4 > 0$ such that 
\begin{align}
	&\E[\Delta_k^2] \leq C_3 \qquad \mbox{ for all }  k \geq 0 ,  \label{eqn:1000-1}\\
	&\E[|\Delta_k|^r] \leq C_4 (a_{k-1}^r + e^{-k}) \qquad \mbox{ for all }  k \geq k_1, \;\; r \in [2,6]. \label{eqn:1000-2}
\end{align}

By Theorem \ref{thm:variance} (ii) and the assumption 
$\sum_{k=2}^\infty (F^{-1}(p_c + 2^{-k}))^2 = \infty$, we have $\var(T(\vzero, \cC_n)) \to \infty$ as $n \to \infty$. 
By \eqref{eqn:1000-1} we have $\sum_{k=0}^{k_1-1} \E[\Delta_{k}^2] \leq C_3 k_1$, thus we can 
throw away the first $k_1$ terms and it is sufficient to prove
\begin{equation*}
	\frac{\sum_{k=k_1}^n \Delta_{k}}{(\sum_{k=k_1}^n \E[\Delta_{k}^2])^{1/2}} \convd N(0,1).
\end{equation*}
From now on we assume $n \geq k_1$. Denote $\sigma_n := (\sum_{k=k_1}^n \E[\Delta_{k}^2])^{1/2}$. Then we have $\lim_{n \to \infty} \sigma_n = \infty$. We verify the conditions in Theorem \ref{thm:mart-clt-mcleish} with $X_{k,n} := \Delta_{k}/\sigma_n$ as follows. Applying Markov's inequality and \eqref{eqn:1000-2} with $r=3$, for any fixed $x > 0$ we have 
\begin{equation}
	\label{eqn:701}
	\pr\left( \max_{k_1  \leq k \leq n } |\Delta_{k}| > x \sigma_n \right) 
	\leq \sum_{k=k_1}^n \pr\left(  |\Delta_{k}| > x \sigma_n \right) 
	\leq \sum_{k=k_1}^n \frac{ C_4 (a_{k-1}^3 + e^{-k})}{x^3 \sigma_n^3} \leq \frac{C_4}{x^3} \cdot \frac{1 +  a_{k_1 - 1} \sum_{k=k_1}^n a_{k-1}^2 }{\sigma_n^3}.
\end{equation}
By Lemma \ref{lem:delta-p-lower-bound}, there is $K_{\ref{lem:delta-p-lower-bound}} > 0$ such that for all $k \geq 2$, we have $\E[\Delta_{k}^2] \geq K_{\ref{lem:delta-p-lower-bound}} [F^{-1}(p_c + 2^{-\delta_0 k})]^2$. Combining this with Lemma \ref{lem:sandwich-series}, there exists a constant $C_5 > 0$ such that for all $n \geq k_1$,
\begin{equation}
	\label{eqn:713}
	\sum_{k=k_1}^n a_{k-1}^2 \leq C_5 \sigma_n^2.
\end{equation}
Combining  \eqref{eqn:713} and \eqref{eqn:701},
\begin{equation*}
	\pr\left( \max_{k_1  \leq k \leq n } |\Delta_{k}| > x \sigma_n \right) \leq \frac{C_4}{x^3} \cdot \frac{1 +  a_{k_1-1} C_5 \sigma_n^2 }{\sigma_n^3}.
\end{equation*}
The above bounds imply conditions (i) and (ii) of Theorem \ref{thm:mart-clt-mcleish}. Next we verify (iii), which states
\begin{equation}
	\label{eqn:716}
		\frac{\sum_{k=k_1}^n (\Delta_{k}^2 - \E[\Delta_{k}^2])}{\sigma_n^2} \convp 0.
\end{equation}
To prove the above convergence it is sufficient to show $\var(\sum_{k=k_1}^n \Delta_{k}^2) = o(\sigma_n^4)$. The following tool, which is \cite[(2.2)]{davydov1968convergence}, gives us a bound on the covariance of terms in this sum.
\begin{lem}[Davydov]
Let $k,m \geq 1$ and let $f,g$ be functions such that $f$ is measurable relative to $\sigma(\Delta_j : j \leq k)$ and $g$ is measurable relative to $\sigma(\Delta_j : j \geq k+m)$. Suppose that $1/p+1/q<1$ and that the moments $\mathbb{E}|f|^p$ and $\mathbb{E}|g|^q$ exist. Then
\[
|\cov(f,g)| \leq 12 \E[|f|^p]^{1/p} \E[|g|^q]^{1/q} \alpha(m)^{1-1/p-1/q}.
\]
\end{lem}
Using this lemma and \eqref{eqn:1000-2} with $r = 6$, for $k \geq k_1$ and $\ell \geq k_1$ we obtain
\begin{align*}
	\cov( \Delta_k^2, \Delta_{\ell}^2 ) 
	\leq& 12 \E[|\Delta_{k}|^{6}]^{1/3}\E[|\Delta_{\ell}|^{6}]^{1/3} [\alpha(|\ell-k|)]^{1/3}  \\
	\leq& 12 C_4^{2/3} (a_{k-1}^2 + e^{-k/3})(a_{\ell-1}^2 + e^{-\ell/3}) [\alpha(|\ell-k|)]^{1/3}.
\end{align*}
Combining the above bound and Lemma \ref{lem:strong-mixing-coefficient}, there exist constants $C_6, C_7>0$ such that for all $\ell, k \geq k_1$,
\begin{equation*}
	\cov(\Delta_{k}^2, \Delta_{\ell}^2)  \leq C_6 e^{-C_7 |k - \ell|} (a_{k-1}^2 + e^{-k/3}) (a_{\ell-1}^2 + e^{-\ell/3}).
\end{equation*}
Therefore we have, as $n \to \infty$,
\begin{align*}
	\var(\sum_{k=k_1}^n \Delta_{k}^2) 
	\leq&  \sum_{k=k_1}^n \sum_{\ell = k_1}^n C_6 e^{-C_7 |k - \ell|} (a_{k-1}^2 + e^{-k/3}) (a_{\ell-1}^2 + e^{-\ell/3})\\
	=& O\bpare{\sum_{k=k_1}^n \sum_{\ell = k_1}^n e^{-C_7 |k - \ell|} a_{k-1}^2 a_{\ell-1}^2} + O\bpare{\sum_{k=k_1}^n a_{k-1}^2} + O(1).
\end{align*}
Since $\sum_{k=k_1}^n a_{k-1}^2 \to \infty$ as $n \to \infty$, one can show that
\begin{equation*}
	\sum_{k=k_1}^n \sum_{\ell = k_1}^n e^{-C_7 |k - \ell|} a_{k-1}^2 a_{\ell-1}^2 = o\bpare{(\sum_{k=k_1}^n a_{k-1}^2)^2} \quad \mbox{ as } n \to \infty.
\end{equation*}
Combining this and \eqref{eqn:713} completes the proof of \eqref{eqn:716}.  \qed

\subsection{Proof of Theorems \ref{thm:variance} and \ref{thm:limit-theorem}.}
\label{sec:proof-variance-clt-box}

In this section, we prove results about $T(\vzero, \partial B(n))$ using Theorems \ref{thm:variance-zero-circuit} and \ref{thm:clt-zero-circuit}.  For any $n \geq 1$, let $q$ be the integer such that $2^{q-1} \leq n < 2^q$. The following lemma controls the error between $T(\vzero, \partial B(n))$ and $T(\vzero, \cC_q)$.

\begin{lem}\label{lem:error-bn-cq}
	Recall $\eta_0$ from \eqref{eqn:def-eta0}. Assume $\eta_0 > 1$. \\
	(i) For any $r \in [1, \eta_0)$, there is $C_0 > 0$ such that for all $n \geq 1$ and $q \geq 1$ such that $2^{q-1} \leq n < 2^q$
	\begin{equation*}
		\E\bbrac{|T(\vzero, \partial B(n)) - T(\vzero,\cC_q)|^r } < C_0.
	\end{equation*}
	(ii) Assume that $\sum_k a_k^{\eta_1} < \infty$ for some $\eta_1 \in [1,\eta_0)$. Then
	\begin{equation*}
		\sum_{q=0}^\infty \; \sup_{ 2^{q-1} \leq n < 2^q }  \E\bbrac{|T(\vzero, \partial B(n)) - T(\vzero,\cC_q)|^{\eta_1} } < \infty.
	\end{equation*}
\end{lem}
\textbf{Proof:} We first prove (i). Observe that, for $1 \leq \ell \leq q $, on the event $\set{m(q-\ell) \geq q-1 > m(q-\ell - 1)}$, $ \partial B(n)$ is sandwiched between $\cC_{q-\ell-1}$ and $\cC_{q}$. Furthermore, for integers $1 \leq \ell \leq q $ and $t \geq 0$, restricted to the event $\set{m(q-\ell) \geq q-1 > m(q-\ell - 1)} \cap \set{m(q) = q+t}$, we have
\begin{equation}
	\label{eqn:955}
	|T(\vzero, \partial B(n)) - T(\vzero,\cC_q)| \leq T(\partial B(2^{q-\ell -1}), \partial B(2^{q+t+1})).
\end{equation}
Then define the events $A_\ell := \set{m(q-\ell) \geq q-1 > m(q-\ell - 1)}$, for $1 \leq \ell \leq q$, and $B_t :=\set{m(q) = q+t}$, for $t \geq 0$. Using \eqref{eqn:955} and the fact that $\cup_{1 \leq \ell \leq q} \cup_{t \geq 0} (A_\ell \cap B_t)$ cover the whole probability space $\Omega$, we have
\begin{align}
	\E\bbrac{|T(\vzero, \partial B(n)) - T(\vzero,\cC_q)|^r } 
	\leq& \sum_{\ell = 1}^q \sum_{t=0}^\infty  \E\bbrac{T^r(\partial B(2^{q-\ell -1}), \partial B(2^{q+t+1})) \ind_{A_{\ell}} \ind_{B_t}}  \nonumber\\
	\leq& \sum_{\ell = 1}^q \sum_{t=0}^\infty  \E\bbrac{T^\eta(\partial B(2^{q-\ell -1}), \partial B(2^{q+t+1}))}^{\frac{r}{\eta}} \pr(A_{\ell})^{\frac{\eta-r}{2\eta}} \pr(B_t)^{\frac{\eta-r}{2\eta}}, \label{eqn:1017}
\end{align}
where the last line uses H\"older's inequality with fixed $\eta \in (r,\eta_0)$. Recall the constant $k_0$ from Lemma \ref{lem:moments-t-k}. Define $b_k := a_k + e^{-k}$ for $k \geq k_0$ and $b_k := b_{k_0}$ for $-1 \leq k < k_0$. By Lemma \ref{lem:483}, there exists $C_1 > 0$ such that for all integers $k \geq -1$ and $r \geq 0$ 
\begin{equation}
	\label{eqn:1021}
	 \E\bbrac{T^\eta(\partial B(2^{k}), \partial B(2^{k+r}))} \leq (C_1 r b_k)^{\eta}.
\end{equation}
By Lemma \ref{lem:delta-p-decomposition}(i), there exists $C_2 > 0$ such that for all integers $1 \leq \ell \leq q$ and $t \geq 0$
\begin{equation}
	\label{eqn:1026}
	\pr(A_{\ell})^{\frac{\eta-1}{2\eta}} \pr(B_t)^{\frac{\eta-1}{2\eta}} \leq e^{-C_2(\ell -1)} e^{-C_2 t}.
\end{equation}
Combining \eqref{eqn:1017}, \eqref{eqn:1021} and \eqref{eqn:1026} we have
\begin{align*}
	\E\bbrac{|T(\vzero, \partial B(n)) - T(\vzero,\cC_q)|^r }  \leq \sum_{\ell = 1}^q \sum_{t=0}^\infty C_1^r (t+\ell + 2)^r b_{q-\ell -1}^r e^{-C_2(\ell -1)} e^{-C_2 t}
	= \sum_{\ell = 1}^q b_{q-\ell -1}^r c_{\ell},
\end{align*}
where $c_\ell := e^{-C_2(\ell -1)}\sum_{t=0}^\infty C_1^r (t+\ell + 2) e^{-C_2 t}$ for $\ell \geq 1$. Write $b_{k} := 0$ for $k \leq -2$ and $c_{\ell} := 0$ for $\ell \leq -1$. Define $\tilde b :=(b_k^r: k \in \bZ)$ and $\tilde c := (c_k: k \in \bZ)$. Then the above bound can be written as $(\tilde b * \tilde c)_{q-1}$, where $\tilde b* \tilde c$ is the convolution of $\tilde b$ and $\tilde c$. Note that $\|\tilde b\|_\infty < \infty$ and $\|\tilde c\|_1 < \infty$. Then (i) follows from Young's inequality, which says $\|\tilde b * \tilde c\|_{\infty} \leq \|\tilde b\|_\infty \|\tilde c\|_1$. 

Next we prove (ii). Replacing $r$ with $\eta_1$ in the above argument, we have $\E\bbrac{|T(\vzero, \partial B(n)) - T(\vzero,\cC_q)|^{\eta_1} } \leq (\tilde b * \tilde c)_{q-1}$. Therefore by Young's inequality,
\begin{equation*}
	\sum_{q=0}^\infty\; \sup_{n: \; 2^{q-1} \leq n < 2^q }  \E\bbrac{|T(\vzero, \partial B(n)) - T(\vzero,\cC_q)|^{\eta_1} }  \leq \|\tilde b * \tilde c\|_1 \leq \|\tilde b\|_1 \|\tilde c\|_1.
\end{equation*}
The assumption $\sum_{k=k_0}^\infty a_k^{\eta_1} < \infty$ implies $\|\tilde b\|_1 < \infty$. Thus the proof of (ii) is completed. \qed\\

Now we are ready to prove our main results about $T(\vzero, \partial B(n))$, beginning with the variance bound.

\textbf{Proof of Theorem \ref{thm:variance}: } For simplicity, denote $s_q := \sum_{k=2}^q [F^{-1}(p_c + 2^{-k})]^2$. For $n \geq 2$, let $q \geq 2$ be the integer such that $2^{q-1} \leq n < 2^{q} - 1$. Denote $X_n := T(\vzero, \partial B(n)) - \E T(\vzero, \partial B(n))$ and $Y_n := T(\vzero, \cC_q) - \E T(\vzero, \cC_q)$. Since $\eta_0 > 2$, we may apply Lemma \ref{lem:error-bn-cq}(i) with $r=2$, there exist a constant $C_0 > 0$ such that for all $n \geq 2$
\begin{equation}
	\label{eqn:1435}
	\|X_n - Y_n\|_2 = \E[|X_n - Y_n|^2]^{1/2} \leq C_0.
\end{equation}
By Theorem \ref{thm:variance-zero-circuit}, there exist $C_1, C_2 > 0$ such that for all $n \geq 2$,
\begin{equation*}
	 C_1 \sqrt{s_q} \leq \| Y_n \|_2 \leq C_2 \sqrt{s_q}.
\end{equation*}
Combining the above two bounds and the triangle inequality, we have
\begin{equation*}
	 ((C_1 \sqrt{s_q} - C_0)\vee 0)^2 \leq \E [X_n^2] \leq (C_2 \sqrt{s_q} + C_0)^2.
\end{equation*}

This suffices to prove the upper bound. For the lower bound, the term $C_1 \sqrt{s_q}-C_0$ may be negative for small $n$, so the proof will be complete once we show that $\var T(\vzero, \partial B(n))>0$ uniformly in $n \geq 1$. The proof of this is standard. Let $S$ be the collection of edges adjacent to $\vzero$ and decompose a configuration $(t_e)$ as $(t_S, t_{S^c})$, where $t_S$ is the collection of edge-weights for edges in $S$ and similarly for $S^c$. Then writing $T = T(\vzero, \partial B(n))$ and $T = T(t_S,t_{S^c})$ to emphasize dependence on the edge-weights, we use Jensen's inequality for
\begin{align*}
\var T &= \int \int \left( T(t_S,t_{S^c}) - \int \int T(t_S', t_{S^c}') ~\text{d}(t_S')~\text{d}(t_{S^c}') \right)^2~\text{d}(t_S)~\text{d}(t_{S^c}) \\
&\geq  \int \left(\int T(t_S,t_{S^c}) ~\text{d}(t_{S^c}) - \int \int T(t_S', t_{S^c}')  ~\text{d}(t_S')~\text{d}(t_{S^c}') \right)^2~\text{d}(t_S) \\
&= \int \left( \int \int [ T(t_S,t_{S^c}) - T(t_S',t_{S^c})]~\text{d}(t_{S^c})~\text{d}(t_S') \right)^2~\text{d}(t_S).
\end{align*}
Since $F(0) = p_c<1$, we can find $\chi>0$ such that $0 < F(\chi) < 1$. Set $A = \{t_e = 0 \text{ for all } e \in S\}$ and $B = \{t_e > \chi \text{ for all } e \in S\}$. We obtain the bound
\[
\var T \geq \int_{t_S \in A} \left( \int \int [T(t_S,t_{S^c}) - T(t_S',t_{S^c})]~\text{d}(t_{S^c})~\text{d}(t_S') \right)^2~\text{d}(t_S).
\]
On the event $\{t_S \in A\}$, the quantity $T(t_S,t_{S^c}) - T(t_S',t_{S^c})$ is nonpositive, so we obtain a bound by restricting to the set $\{t_S' \in B\}$, on which $T(t_S,t_{S^c}) - T(t_S',t_{S^c}) \leq -\chi$:
\begin{align*}
\var T &\geq \int_{t_S \in A} \left( \int_{t_S' \in B} \int [T(t_S,t_{S^c}) - T(t_S',t_{S^c})]~\text{d}(t_{S^c})~\text{d}(t_S') \right)^2~\text{d}(t_S) \\
&\geq (\chi \mathbb{P}(t_S \in B))^2 \mathbb{P}(t_S \in A) > 0.
\end{align*}
This bound does not depend on $n$, so we are done.
 \qed\\


Last, we deduce limit theorems for $T(\vzero, \partial B(n))$.

\textbf{Proof of Theorem \ref{thm:limit-theorem}: } First we prove (i).
 Suppose $\sum_{k = 2 }^\infty [F^{-1}(p_c + 2^{-k})]^2 < \infty$. Then $\sum_{k=k_0}^\infty a_k^2 < \infty$, where $k_0$ is defined in Lemma \ref{lem:moments-t-k}. Also note that $T(\vzero,\cC_q) - \E T(\vzero, \cC_q) = \sum_{k=0}^q \Delta_k$ and $\Delta_k$, for $k \geq 0$, does not depend on $n$ or $q$. By Theorem \ref{thm:variance-zero-circuit}(ii), we have $\sum_{k=1}^\infty \E \Delta_k^2 < \infty$. Then by the Martingale Convergence Theorem, there exists a random variable $Z$ with  $\E Z = 0$ and $\E Z^2 < \infty$  such that as $q \to \infty$
\begin{equation}
	\label{eqn:1071}
		T(\vzero,\cC_q) - \E T(\vzero,\cC_q) \to Z \quad \mbox{ a.s. and in $L^2$}.
\end{equation}
Applying Lemma \ref{lem:error-bn-cq}(ii) with $\eta_1 = 2$ and taking $n_q = 2^{q-1}$ or $2^{q}-1$, for $q \geq 0$, we have
\begin{equation*}
		\sum_{q=0}^\infty \E\bbrac{|T(\vzero, \partial B(n_q)) - T(\vzero,\cC_q)|^2 } < \infty.
\end{equation*}
Therefore by Borel-Cantelli and \eqref{eqn:1071}, as $q \to \infty$,
\begin{equation}
	\label{eqn:1080-1}
		T(\vzero,\partial B(n_q)) - \E T(\vzero,\partial B(n_q)) \to Z \quad \mbox{ a.s. and in $L^2$}.
\end{equation}
Note that for all $n,q$ such that $2^{q-1} \leq n < 2^q$ we have
\begin{equation*}
	|T(\vzero, \partial B(n)) - T(\vzero, \cC_q)|  \leq \max\set{ |T(\vzero, \partial B(2^{q-1})) - T(\vzero, \cC_q)|, |T(\vzero, \partial B(2^q-1)) - T(\vzero, \cC_q)| }.
\end{equation*}
Combining the above observation and \eqref{eqn:1080-1} completes the proof of Theorem \ref{thm:limit-theorem}(i).

Next we prove (ii). Suppose $\sum_{k = 2 }^\infty [F^{-1}(p_c + 2^{-k})]^2 = \infty$. Define $\sigma_n := \var(T(\vzero,\cC_q))^{1/2}$ where $q \in \bN$ is such that $2^{q-1} \leq n < 2^q$. Define $\gamma_n := \var(T(\vzero,  \partial B(n)))$. By Theorem \ref{thm:variance-zero-circuit}(ii) we have $\lim_{n \to \infty} \sigma_n = \infty$.  By Lemma \ref{lem:error-bn-cq}(i) with $r=2$, there is $C_0 >0$ such that for all $n \geq 2$
\begin{equation}
	\label{eqn:1488}
	|\sigma_n - \gamma_n| \leq C_0.
\end{equation}
Furthermore, there is $C_1 > 0$ such that for all $n \ge 2$
\begin{equation}
	\label{eqn:1492}
	\E\bbrac{|T(\vzero, \partial B(n)) - T(\vzero,\cC_q)|} \leq C_1.
\end{equation}
Then Theorem \ref{thm:limit-theorem}(ii) is a consequence of Theorem \ref{thm:clt-zero-circuit}, \eqref{eqn:1492}, \eqref{eqn:1488} and the fact that $\lim_{n \to \infty} \sigma_n = \infty$. \qed\\

\subsection{Limit theorems for point-to-point times}
\label{sec: limit-0-to-x}



In this section we extend the limit theorems and variance estimates from the last section to point-to-point passage times. 

\begin{cor}
	\label{cor:mean-p2p}
	(i) Assume that $\eta_0 > 1$. There exists $C_1 = C_1(F) > 0$ such that
	\begin{equation*}
		\E T(\vzero, x) \leq C_1 \sum_{k = 2}^q F^{-1}(p_c + 2^{-k}) \qquad \mbox{ for } \;\; x \in B(2^{q+1})\setminus B(2^{q}) \;\; \mbox{ and } \;\; q \geq 2.
	\end{equation*}
	(ii) There exists $C_2 = C_2(F) > 0$ such that
	\begin{equation*}
		\E T(\vzero, x) \geq C_2 \sum_{k = 2}^q F^{-1}(p_c + 2^{-k}) \qquad \mbox{ for } \;\; x \in B(2^{q+1})\setminus B(2^{q}) \;\; \mbox{ and } \;\; q \geq 2.
	\end{equation*}
\end{cor}

\begin{cor} 
	\label{cor:variance-p2p}
	Assume that $\eta_0 > 2$.\\
	(i)  There exists $C_3 = C_3(F) > 0$ such that
	\begin{equation*}
		\var(T(\vzero, x)) \leq C_3 \sum_{k = 2}^q [F^{-1}(p_c + 2^{-k})]^2 \qquad \mbox{ for } \;\; x \in B(2^{q+1})\setminus B(2^{q}) \;\; \mbox{ and } \;\; q \geq 2.
	\end{equation*}
	(ii) There exists $C_4 = C_4(F) > 0$ such that
	\begin{equation*}
		\var(T(\vzero, x)) \geq C_4 \sum_{k = 2}^q [F^{-1}(p_c + 2^{-k})]^2 \qquad \mbox{ for } \;\; x \in B(2^{q+1})\setminus B(2^{q}) \;\; \mbox{ and } \;\; q \geq 2.
	\end{equation*}
\end{cor}

\begin{cor} 
	\label{cor:limit-theorem-p2p}
	Assume that $\eta_0 > 2$ and $\sum_{k=2}^\infty F^{-1}(p_c + 2^{-k}) = \infty$.\\
	(i) If $\sum_{k=2}^\infty [F^{-1}(p_c + 2^{-k})]^2 < \infty$, then there exists a random variable $\tilde Z$ with $\E \tilde Z = 0$ and $\E \tilde Z^2 < \infty$ such that
	\begin{equation*}
		T(\vzero, x) - \E T(\vzero, x) \convd \tilde Z   \qquad \mbox{ as } \;\; \|x\|_{\infty} \to \infty.
	\end{equation*}
	Here $\tilde Z$ has the same distribution as the sum of two independent copies of $Z$, which is defined in Theorem \ref{thm:limit-theorem}.\\
	(ii) If $\sum_{k=2}^\infty [F^{-1}(p_c + 2^{-k})]^2 = \infty$, then
	\begin{equation*}
		\frac{T(\vzero, x) - \E T(\vzero, x)}{\var(T(\vzero, x))^{1/2}} \convd N(0,1) \qquad \mbox{ as } \;\; \|x\|_{\infty} \to \infty. 
	\end{equation*}
	In particular, letting $q=q(x)$ be the integer such that $2^{q} < \|x\|_\infty \leq 2^{q+1}$, we have
	\begin{equation*}
		\frac{\var(T(\vzero,x))}{ \var(T(\vzero,\partial B(2^{q(x)})))} \to 2 \qquad \mbox{ as } \|x\|_\infty \to \infty.
	\end{equation*}
\end{cor}

\begin{rem}
	In comparison to the $L^2$ and a.s. convergence in Theorem \ref{thm:limit-theorem}(ii), one would only expect convergence in distribution in Corollary~ \ref{cor:limit-theorem-p2p} (i). This can be explained by the following fact: $T(\vzero,x)$ heavily depends on the edge-weights near the point $x$, which tends to infinity. As $x$ changes, the edge weights near it only share the same distribution.
\end{rem}

Now we describe the main construction that is used in the proof of the above three corollaries. This construction was introduced in \cite{kesten1997central}. Suppose $x \in B(2^{q+1})\setminus B(2^{q})$. Then the two boxes $B(\vzero,2^{q-1})$ and $B(x,2^{q-1})$ are disjoint, and therefore the two random variables $T(\vzero, \partial B(\vzero,2^{q-1}))$ and $T(x, \partial B(x,2^{q-1}))$ are independent and identically distributed. Define
\begin{equation*}
	Y(x) := T(\vzero, \partial B(\vzero,2^{q-1})) + T(x, \partial B(x,2^{q-1})) \qquad \mbox{ for } x \in B(2^{q+1})\setminus B(2^q) \mbox{ and } q \geq 2.
\end{equation*}
Then $T(\vzero, x) \geq Y(x)$. Note that the statements in the above three corollaries, with $T(\vzero,x)$ replaced by $Y(x)$, are all immediate consequences of Theorems \ref{thm:mean}, \ref{thm:variance} and \ref{thm:limit-theorem}. Then we only need to control the error between $T(\vzero,x)$ and $Y(x)$. To bound $T(\vzero,x)$ from above, recall the definition of $\cC_{q+2}$ from \eqref{eqn:def-c-n}. One can construct a path between $\vzero$ and $x$ by concantenating a geodesic from $\vzero$ to $\cC_{q+2}$, a $p_c$-open path along $\cC_{q+2}$, and a geodesic from $\cC_{q+2}$ to $x$. Thus $T(\vzero, x)$ can be bounded above by $T(\vzero, \cC_{q+2}) + T(x, \cC_{q+2})$. This implies
\begin{equation}
	\label{eqn:1940}
	|T(\vzero, x) - Y(x)| \leq |T(\vzero, \cC_{q+2}) - T(\vzero, \partial B(\vzero,2^{q-1}))|+ |T(x, \cC_{q+2}) - T(x, \partial B(x,2^{q-1}))|.
\end{equation}
The first term in the above bound can be controlled by Lemma \ref{lem:error-bn-cq} and the second term can be controlled by the following lemma, which is also analogous to Lemma \ref{lem:error-bn-cq}.
\begin{lem}
	\label{lem:1945}
	Recall $\eta_0$ from \eqref{eqn:def-eta0}. Assume $\eta_0 > 1$. \\
	(i) For any $r \in [1, \eta_0)$, there is $C_0 > 0$ such that for all $q \geq 0$ and $x \in B(2^{q+1})\setminus B(2^q)$
	\begin{equation*}
		\E\bbrac{|T(x, \cC_{q+2}) - T(x, \partial B(x,2^{q-1}))|^r } < C_0.
	\end{equation*}
	(ii) Assume that $\sum_k a_k^{\eta_1} < \infty$ for some $\eta_1 \in [1,\eta_0)$. Then
	\begin{equation*}
		\sum_{q=0}^\infty \; \sup_{ x \in B(2^{q+1})\setminus B(2^q) }  \E\bbrac{|T(x, \cC_{q+2}) - T(x, \partial B(x,2^{q-1}))|^{\eta_1} } < \infty.
	\end{equation*}
\end{lem}
The proof of the above lemma is similar to the one of Lemma \ref{lem:error-bn-cq}, and therefore is omitted.

\textbf{Proof of Corollary \ref{cor:mean-p2p}:} By Lemma \ref{lem:1945}(i), Lemma \ref{lem:error-bn-cq}(i) and \eqref{eqn:1940}, there exists a constant $C_0 > 0$ such that $\E |T(\vzero, x) - Y(x)| \leq C_0$ for all $x$. This proves (i). Combining the lower bound $T(\vzero, x) \geq Y(x)$ and Theorem \ref{thm:mean}(ii) proves (ii). \qed \\

\textbf{Proof of Corollary \ref{cor:variance-p2p}:} When $\eta_0 > 2$, by  Lemma \ref{lem:1945}(i), Lemma \ref{lem:error-bn-cq}(i) and \eqref{eqn:1940}, there exists a constant $C_0 > 0$ such that $\E |T(\vzero, x) - Y(x)|^2 \leq C_0$ for all $x$. Then the rest of the proof is similar to the proof of Theorem \ref{thm:variance}. \qed \\

\textbf{Proof of Corollary \ref{cor:limit-theorem-p2p}:} To show (i), since $\sum_{k} a_k^2 < \infty$ and $\eta_0 > 2$, then by Lemma \ref{lem:1945}(ii) we have
$\E|T(x, \cC_{q+2}) - T(x, \partial B(x,2^{q-1}))|^2 \to 0$ as $\|x\|_{\infty} \to \infty$. Then by \eqref{eqn:1940} we have
\begin{equation*}
	T(\vzero,x) - Y(x) \to 0 \qquad \mbox{ in $L^2$ as } \|x\|_\infty \to \infty.
\end{equation*} 
By Theorem \ref{thm:limit-theorem}(i) and the independence of $T(\vzero, \partial B(2^{q-1}))$ and $T(x, \partial B(x, 2^{q-1}))$, we have 
\begin{equation*}
	Y(x) \convd Z + Z', \qquad \mbox{ as }  \|x\|_\infty \to \infty,
\end{equation*}
where $Z'$ is another independent copy of $Z$ as in Theorem \ref{thm:limit-theorem}(i). Combining these proves (i). The proof of (ii) is similar to that of Theorem \ref{thm:limit-theorem}(ii). \qed\\

\bigskip
\noindent
{\bf Acknowledgements.} The research of M. D. is supported by NSF grant DMS-1419230.

%
%



\bibliographystyle{plain}
\bibliography{first_passage_percolation}

\end{document}